\documentclass[11pt,a4paper]{article}
\usepackage[centertags]{amsmath}
\usepackage{amsfonts}
\usepackage{amssymb}
\usepackage{amsthm}
\usepackage{newlfont}

\newlength{\defbaselineskip}
\theoremstyle{plain}
\newtheorem{thm}{Theorem}[section]
\newtheorem{lem}{Lemma}[section]
\newtheorem{prop}{Proposition}[section]
\newtheorem{defn}{Definition}[section]
\newtheorem{cor}{Corollary}[section]
\newtheorem{rem}{Remark}[section]
\numberwithin{equation}{section}

\newcommand{\A}{\alpha}
\newcommand{\B}{\beta}
\newcommand{\C}{\gamma}

\newcommand{\M}{\overline{M}}
\newcommand{\MM}{\widetilde{M}}
\newcommand{\g}{\overline{g}}
\newcommand{\G}{\widetilde{g}}

\newcommand{\n}{\nabla}
\newcommand{\p}{\partial}

\newcommand{\PP}{\widetilde{P}}

\title{Induced structures on submanifolds in almost product Riemannian
manifolds}
\author{Cristina-Elena Hre\c{t}canu}
\date{}

\usefont{T1}{cmss}{m}{ol}

\begin{document}
 \maketitle

\begin{abstract}
We give some fundamental properties of the induced structures on
submanifolds immersed in almost product or locally product
Riemannian manifolds. We study the induced structure by the
composition of two isometric immersions on submanifolds in an
almost product Riemannian manifold. We give an effective
construction for some induced structures on submanifolds of
codimension 1 or 2 in Euclidean space.
\end{abstract}

\section*{Introduction}
The geometry of submanifolds with induced structures in Riemannian
manifolds was widely studied by many geometers, such as K. Yano
and M. Kon (\cite{KY}, \cite{Yano1}, \cite{Yano2}, \cite{Yano2}).
An investigation of the properties of the almost product or
locally product Riemannian manifolds has been made by M. Okumura
(\cite{Okumura}), T. Adati and T. Miyazawa (\cite{Adati1},
\cite{Adati2}), M. Anastasiei(\cite{AM1}), G. Piti\c{s}
(\cite{Pitis1}), X. Senlin and N. Yilong (\cite{Senlin}), A. G.
Walker (\cite{Walker}), M. At\c{c}eken, S. Kele\c{s} and B.
\c{S}ahin (\cite{Atceken}, \cite{Sahin}), etc. Also, the
properties of the almost r-paracontact structures were studied by
A. Bucki and A. Miernovski (\cite{Bucki1}, \cite{Bucki2}), T.
Adati and T. Miyazawa (\cite{Adati3}), S. Ianu\c{s} and I.Mihai
(\cite{IS7}), J. Nikic (\cite{Nikic}), etc.

The purpose of this paper is to give some properties of the
submanifolds with a
$(P,g,\varepsilon\xi_{\A},u_{\A},(a_{\A\B})_{r})$ structure
induced by a $\PP$ structure defined on a Riemannian manifold
$(\MM,\G)$, with $\PP^{2}=\varepsilon I$ and the compatibility
equality (1.2) between $\G$ and $\PP$ (where $I$ is the identity
on $\MM$ and $\varepsilon = \pm 1$). The
$(P,g,\varepsilon\xi_{\A},u_{\A},(a_{\A\B})_{r})$ structure is
determined by an (1,1)-tensor field P on M, tangent vector fields
$\xi_{\A}$ on $M$, 1-forms $u_{\A}$  on $M$ and a $r \times r$
matrix $(a_{\A\B})_{r}$, where its entries $a_{\A\B}$ are real
functions on $M$ ($\A, \B \in \{1,...,r\}$). Particularly, for $
\varepsilon =1$, the $\PP$ structure on $(\MM,\G)$ becomes an
almost product structure.

This paper is organized as follows: in section 1 we construct the
structure $(P,g,\varepsilon \xi_{\A},u_{\A},(a_{\A\B})_{r})$,
induced on a submanifold in the Riemannian manifold $(\MM,\G,\PP)$
with the conditions (1.1) and (1.2), in the same manner like in
\cite{Adati2}.

In section 2, we give the fundamental formulae for
$(P,g,\varepsilon\xi_{\A},u_{\A},(a_{\A\B})_{r})$ induced
structure on a submanifold in the Riemannian manifold
$(\MM,\G,\PP)$, with $\widetilde{\n}\PP=0$.

In sections 3 and 4, we shall investigate the necessary and
sufficient conditions for a submanifold with
$(P,g,\xi_{\A},u_{\A},(a_{\A\B})_{r})$ induced structure, immersed
in a locally Riemannian product manifold to be normal (relative to
the commutativity of the endomorphism P and the Weingarten
operators $A_{\A}$ on M) and show further properties of this kind
of submanifold.

In section 5, we prove that the composition of two isometric
immersions
$M\hookrightarrow\overline{M}\hookrightarrow\widetilde{M}$ induces
on a submanifold M of codimension 2 in an almost product
Riemannian manifold $(\widetilde{M},\widetilde{g},\widetilde{P})$
a $(P,g,u_{1},u_{2},\xi_{1},\xi_{2}^{\bot},(a_{\alpha\beta}))$
structure determined by a
$(\overline{P},\overline{g},u_{2},\xi_{2},a_{22})$ structure on
$\overline{M}$ (induced by $(\PP,\G)$) and a
$(P,g,u_{1},\xi_{1},a_{11})$ structure on M (induced by
$(\overline{P},\overline{g},u_{2},\xi_{2},a_{22})$), where
$a_{12}=a_{21}=g(\xi_{2}^{\bot},N_{1})$ and $N_{1}$ is a unit
normal vector field on M.

In section 6, we show some properties of
$(P,g,\xi_{\A},u_{\A},(a_{\A\B})_{r})$ induced structures on the
submanifold of codimension 1 or 2 in an  almost product Riemannian
manifold $(\widetilde{M},\widetilde{g},\widetilde{P})$.

In section 7, we give an effective construction for some
$(P,g,\xi_{\A},u_{\A},(a_{\A\B})_{r})$ induced structures on
hyperspheres or submanifolds of codimension 2 in Euclidean  space.

I would like to express here my sincere gratitude to Professor Dr.
Mihai Anastasiei (from "Al.I.Cuza" University - Ia\c{s}i, Romania)
who gave me many suggestions to improve the first draft of this
paper.

\section{\bf $(P,g,\varepsilon
\xi_{\A},u_{\A},(a_{\A\B})_{r})$ induced structure on submanifolds
in Riemannian manifold}

\normalfont Let $(\MM,\G)$ be a Riemannian manifold, equipped with
a Riemannian metric tensor $\G$ and a (1,1) tensor field $\PP$
such that
\[\PP^{2} = \varepsilon I, \:\: \varepsilon=\pm 1 \leqno(1.1)\]
where $I$ is the identity on $\MM$.  We suppose that $\G$ and
$\PP$ are compatible in the sense that for each $U, V \in
\chi(\MM)$ we have that
\[\G(\PP U, \PP V)=\G(U, V), \leqno(1.2) \]
which is equivalent with
\[ \G(\PP U, V)=\varepsilon\G(U, \PP V), \:(\forall) U,V \in \chi(\MM) \leqno(1.3)\]
for each $U, V \in \chi(\MM)$, where $\chi(M)$ is the Lie algebra
of the vector fields on $\MM$.

For $\varepsilon =1$, $\PP$ is an almost product structure and the
Riemannian manifold $(\MM,\G,\PP)$, with the compatibility
relation (1.2), becomes an almost product Riemannian manifold.

Let M be an n-dimensional submanifold of codimension r in a
Riemannian manifold $(\MM,\G,\PP)$ which satisfied the relations
(1.1) and (1.2).

We make the following notations throughout all of this paper: X,
Y, Z, ... are tangential vector fields on M. We denote the tangent
space of M at $x \in M$ by $T_{x}(M)$ and the normal space of M in
x by $T_{x}^{\bot}(M)$. Let $(N_{1},...,N_{r}):=(N_{\A})$ be an
orthonormal basis in $T_{x}^{\bot}(M)$, for every $x \in M$. In
the following statements, the indices range is fixed in this way:
$ \A, \B, \C ...\in \{ 1,...,r\}$. We shall use the Einstein
convention for summation.

The decomposition of the vector fields $\PP X$ and $\PP N_{\A}$
respectively, in the tangential and normal components of M has the
form:
\[ \PP X = PX + \sum_{\A}u_{\A}(X)N_{\A},  \leqno(1.4)\]
for any $X \in \chi(M)$ and
\[ \PP N_{\A}= \varepsilon \xi_{\A} + \sum_{\B}a_{\A \B}
N_{\B}, \quad(\varepsilon = \pm 1)\leqno(1.5)\] where P is a (1,1)
tensor field on M, $\xi_{\A}$ are tangent vector fields on M,
$u_{\A}$ are 1-forms on M and $(a_{\A \B})_{r}$ is an $r \times
r$-matrix and its entries $a_{\A\B}$ are real functions on M. If
$\varepsilon = 1$, the formulae in the next theorem were
demonstrated by T.Adati (in \cite{Adati2}).

\thm Let  M be an n-dimensional submanifold of codimension r in a
Riemannian manifold $(\MM, \G)$, equipped by an (1,1)-tensor field
$\PP$, such that $\G$ and $\PP$ verify the conditions (1.1) and
(1.2). The $(\MM, \G, \PP)$ structure induces on the submanifold M
a $(P,g,u_{\A},\varepsilon\xi_{\A},(a_{\A\B})_{r})$ Riemannian
structure which verifies the following properties:
\[ \begin{cases}
    (i) \quad P^{2}X = \varepsilon (X -\sum_{\A}u_{\A}(X)\xi_{\A}),  \\
    (ii) \quad u_{\A}(PX)= -\sum_{\B}a_{\B\A}u_{\B}(X),\\
    (iii) \quad a_{\A\B} = \varepsilon a_{\B \A}, \\
    (iv) \quad u_{\A}(\xi_{\B})=g(\xi_{\A},\xi_{\B})=
    \delta_{\A \B}-\varepsilon \sum_{\gamma}a_{\A \gamma}a_{\gamma\B},\\
    (v) \quad P\xi_{\A}=-\sum_{\B}a_{\A \B}\xi_{\B}
\end{cases} \leqno(1.6)\] and
\[ \begin{cases}
    (i) \quad u_{\A}(X)=g(X,\xi_{\A}),\\
    (ii) \quad g(PX,Y)=\varepsilon g(X,PY),\\
    (iii) \quad g(PX,PY)= g(X,Y)-\sum_{\A}u_{\A}(X)u_{\A}(Y),
\end{cases} \leqno(1.7) \] for any $X,Y \in \chi(M)$.\\
\normalfont\textbf{Proof:} Applying $\PP$ in the equality (1.4) we
obtain that
\[\PP^{2}X=\PP(PX)+\sum_{\A=1}^{r}u_{\A}(X)\PP(N_{\A}),\]
for every $X \in \chi(M)$.

From (1.1), (1.4) and (1.5) we have
\[ \varepsilon X = P^{2}X +\varepsilon \sum_{\A}u_{\A}(X)\xi_{\A}+
\sum_{\B}(u_{\B}(PX)+\sum_{\A}u_{\A}(X)a_{\A \B})N_{\B},\] for
every $X \in \chi(M)$, and from this, it results (i) and (ii) from
the relations (1.6). Furthermore, the equality (i) can be written
in the following form
\[ P^{2} = \varepsilon (I -\sum_{\A}u_{\A}\otimes \xi_{\A}) \leqno(i)'\]

Applying the equality (1.3) to the normal vector fields $N_{\A}$
and $N_{\B}$, respectively and using the equality (1.5) follows
that \[ \G(\varepsilon \xi_{\A}+\sum_{\gamma =1}a_{\A
\gamma}N_{\gamma}, N_{\B}) = \varepsilon \G(N_{\A}, \varepsilon
\xi_{\B}+\sum_{\gamma=1}a_{\B \gamma}N_{\gamma})\] and from this
we obtain the equality (iii) from (1.6).

From $\PP^{2}N_{\A}=\varepsilon N_{\A}$, using the relations (1.4)
and (1.5) we obtain
\[\varepsilon N_{\A}=\PP^{2}N_{\A}=
\PP(\varepsilon \xi_{\A} +\sum_{\B}a_{\A\B}N_{\B})=\varepsilon \PP
\xi_{\A} +\sum_{\B}a_{\A\B}\PP N_{\B}=\]
\[=\varepsilon (P\xi_{\A}+\sum_{\B}a_{\A\B}\xi_{\B})
+\sum_{\B}[\varepsilon u_{\B}(\xi_{\A}) +
\sum_{\gamma}a_{\A\gamma}a_{\gamma\B}]N_{\B}\] so
\[\varepsilon N_{\A}=\varepsilon
(P\xi_{\A}+\sum_{\B}a_{\A\B}\xi_{\B}) +\sum_{\B}[\varepsilon
u_{\B}(\xi_{\A}) + \sum_{\gamma}a_{\A\gamma}a_{\gamma\B}]N_{\B}\]
Identifying  the tangential components from the last equality we
obtain (v) from (1.6) and identifying the normal components from
the last equality we obtain (iv) from (1.6). Applying the equality
(1.3) to the vector fields $X$ and $N_{\A}$ respectively, we
obtain
\[\G(\PP X,N_{\A})=\varepsilon \G(X,\PP N_{\A})\] and from this it follows that
\[\G(PX +\sum_{\B}u_{\B}(X)N_{\B},N_{\A})=
\varepsilon \G(X, \varepsilon \xi_{\A}+\sum_{\B}a_{\A\B}N_{\B})\]
for any tangent vector fields X on M, so we obtain the equality
(i)(1.7).

Applying the relations (1.3) and (1.4) to the tangential vector
fields X and Y on M we obtain
\[ g(PX,Y)=g(\PP X, Y)=g(\PP^{2} X, \PP Y)=\varepsilon g(X,\PP Y)=
\varepsilon g(X,PY)\] and from this we have (ii) from (1.7).\\

Replacing Y by PY in the equality (ii) from (1.7) and using the
equality (i) from (1.6) we obtain
\[ g(PX,PY)= \varepsilon g(X,P^{2}Y)=\varepsilon^{2} (g(X,Y)-
\sum_{\A=1}u_{\A}(Y)g(X,\xi_{\A})),\] for any $X,Y \in \chi(M)$.
From $\varepsilon^{2}=1$ and using (i) from (1.7), we obtain the
equality (iii) from (7). $ \square $

\rem \normalfont Particularly, for $\varepsilon = -1$ and $a=0$
(if we omit the metric g), we obtain a $(P,u_{\A},-\xi_{\A})$
induced structure, which has the following properties:
\[ \begin{cases}
    (i)\quad P^{2}X =-X+\sum_{\A}u_{\A}(X)\xi_{\A}, \: (\forall)X \in \chi(M)\\
    (ii)\quad u_{\A}\circ P=0, \\
    (iii) \quad u_{\A}(\xi_{\B})=\delta_{\A\B}, \\
    (iv)\quad P\xi_{\A} =0.
\end{cases} \leqno(1.8) \]
Applying P in the equality (i) from (1.8) we obtain
\[ P^{3}+P=0 \leqno(1.9)\]

Therefore, P is an $f(3,1)$-structure and the manifold M endowed
with a $(P,u_{\A},-\xi_{\A})$ structure is a framed manifold.
Besides, we can call this structure an almost r-contact structure.
If r=1, then we call this structure an almost contact structure.

\defn For $\varepsilon = -1$ and $a=0$, the $(P,g,u_{\A},-\xi_{\A})$ induced structure by
 $\PP$ from $(\MM, \G)$, with the properties (1.7) and
(1.8) is called an $f(3,1)$ Riemannian structure.

\rem \normalfont Particularly, for $\varepsilon = 1$ and $a=0$ the
$(P,u_{\A},\xi_{\A})$ induced structure on M has the following
properties:
\[ \begin{cases}
    (i) \quad P^{2}X =X-\sum_{\A}u_{\A}(X)\xi_{\A}, \: (\forall)X \in \chi(M)\\
    (ii) \quad u_{\A}\circ P=0,\\
    (iii)\quad u_{\A}(\xi_{\B})=\delta_{\A\B}, \\
    (iv) \quad P\xi_{\A} =0.
\end{cases} \leqno(1.10) \]
Applying P in the equality (i) from (1.10) we obtain
\[ P^{3}-P=0 \leqno(1.11)\]
Therefore, P is an $f(3,-1)$ structure on M, and
$(P,u_{\A},\xi_{\A})$ induced structure on M, with the properties
(1.10) is called an $f(3,-1)$ framed structure. Besides, this kind
of structure is called an almost r-paracontact structure. For r=1
we obtain an almost paracontact structure.

\defn For $\varepsilon =1$ and $a=0$, the $(P,g,u_{\A},\xi_{\A})$ induced structure on M, with the
properties (1.10) and (1.7) is an $f(3,-1)$ Riemannian structure.

\rem \normalfont Therefore, the $(P,g,u_{\A},\varepsilon
\xi_{\A},(a_{\A\B})_{r})$ induced structure on M by the $\PP$
structure on $(\MM,\G)$ (which verifies the relations (1.1) and
(1.2)) is a generalization of an almost r-contact Riemannian
structure and an almost r-paracontact Riemannian structure,
respectively.

Thus, we have the following situations:\\
  (I): For $\varepsilon =-1$ and $a=0$, we obtain an
 $f(3,1)$ structure and the structure $(P,g,u_{\A},-\xi_{\A})$
becomes an almost r-contact Riemannian structure;\\
 (II): For $\varepsilon =1$ and $a=0$, we obtain an
 $f(3,-1)$ structure and the structure $(P,g,u_{\A},\xi_{\A})$
becomes an almost r-paracontact Riemannian structure.

\defn A $(P,g,u_{\A},\varepsilon \xi_{\A},(a_{\A \B})_{r})$ structure on
a submanifold M of codimension r in a $(\MM,\G,\PP)$ Riemannian
manifold with the proprieties (1.1) and (1.2), which verifies the
properties (1.6) and (1.7)
 is called an $(a,\varepsilon)f$ Riemannian structure.

\rem \normalfont The case of the $(a,-1) f$ Riemannian structure
was studied by K. Yano and M. Okumura (in \cite{Yano1} and
\cite{Yano2}).

In the following issue, we suppose that $\varepsilon =1$.
Therefore, the $\PP$ structure on the Riemannian manifold
$(\MM,\G)$, which satisfies the relations (1.1) and (1.2), is an
almost product structure.

\rem \normalfont If we suppose that $\xi_{1},...,\xi_{r}$ are
linearly independent tangent vector fields on M, it follows that
the 1-forms $u_{1},...,u_{r}$ are linearly independent, too.  The
equality
 \[\sum_{\A=1}^{r}\lambda^{\A}u_{\A}(X)=0\]
 is equivalent with
 \[0=\sum_{\A}\lambda^{\A}g(X,\xi_{\A})=g(X,\sum_{\A}\lambda^{\A}\xi_{\A}), \: (\forall) X \in \chi(M) \]
 thus \[ \sum_{\A=1}^{r}\lambda^{\A}\xi_{\A}=0 \Rightarrow \lambda^{\A}=0 \]
 and from this we obtain that $u_{1},...,u_{r}$ are linearly
 independent on M.

\rem \normalfont We denote by
 \[D_{x}=\{X_{x} \in T_{x}M \: : \: u_{\A}(X_{x})=0\}, \leqno(1.12)\]
for any $\A \in \{1,...,r\}$. We remark that $D_{x}$ is an
$(n-r)$-dimensional subspace in $T_{x}M$ and the function
\[D : x \mapsto D_{x}, \: (\forall) x \in M \leqno(1.13)\] is a
distribution locally defined on M. If $X \in D$, from (1.6)(ii) we
have
\[u_{\A}(PX)=-\sum_{\B}a_{\B\A}u_{\B}(X)=0, \leqno(1.14)\]
for any $ X \in D $, then $PX \in D$. Therefore $D$ is an
invariant distribution with respect to P.

If $D_{x}^{\bot}$ is an orthogonal supplement of $D_{x}$ in
$T_{x}M$, then we obtain the distribution $D^{\bot}: x \mapsto
D_{x}^{\bot}$. Furthermore, we have the decomposition
\[T_{x}M=D_{x} \oplus D_{x}^{\bot},  \leqno(1.15)\] in any point $ x \in
M$. From (1.7)(i) it follows that the vector fields $\xi_{\A}\neq
0$ are orthogonal on $D_{x}$ and $\xi_{\A} \in D_{x}^{\bot}$.
Thus, if $\xi_{\A}\neq 0$ for any $\A \in \{1,...,r\}$, then
$D_{x}^{\bot}$ is generated by $\xi_{1},...,\xi_{r}$ and
$D_{x}^{\bot}$ is r-dimensional in $T_{x}M$.

From (1.6)(v) we remark that the space $D_{x}$ is P-invariant and
P satisfies
\[P^{2}X=X, \leqno(1.16)\]
and
\[ g(PX,PY)=g(X,Y),  \leqno(1.17)\]
for all $X,Y \in D$. Thus P is an almost product Riemannian
structure on D. Furthermore, $rank(P)=n-r$ on D and its
eigenvalues are 1 and -1.

\rem \normalfont Let $\{N_{1},...,N_{r}\}$ and
$\{N'_{1},...,N'_{r}\}$ be two orthonormal basis on a normal space
$T_{x}^{\bot}M$. The decomposition of $N'_{\A}$ in the basis
$\{N_{1},...,N_{r}\}$ is the following
\[ N'_{\A}=\sum_{\gamma=1}^{r}k^{\gamma}_{\A}N_{\gamma},
\leqno(1.18)\]for any $\A \in \{1,...,r\}$, where
$(k^{\gamma}_{\A})$ is an $r \times r$ orthogonal matrix and we
have (from \cite{Adati2}):
\[u'_{\A}=\sum_{\gamma}k^{\gamma}_{\A}u_{\gamma} \leqno(1.19)\]
\[\xi'_{\A}=\sum_{\gamma}k^{\gamma}_{\A}\xi_{\gamma} \leqno(1.20)\]
and
\[a'_{\A\B}=\sum_{\gamma}k^{\gamma}_{\A}a_{\gamma \delta}k^{\delta}_{\B} \leqno(1.21)\]
 From (1.20) we have that if $\xi_{1},...,\xi_{r}$ are linearly independent vector fields, then $\xi'_{1},...,\xi_{r}$
 are also linearly independent.


\section{\bf The fundamental equations of submanifolds with $(P,g,\xi_{\A},u_{\A},(a_{\A\B})_{r})$ structures}

\normalfont In this section, we suppose that the Riemannian
manifolds $(\MM,\G)$ are endowed with a (1,1) tensor field $\PP$
on $\MM$, which verifies the equalities (1.1),(1.2), and the
structure $\PP$ is parallel with respect to the Levi-Civita
connection $\widetilde{\n}$ of $\G$. Let M be an n-dimensional
Riemannian submanifold of codimension r, isometric immersed in
$\MM$ and $(P,g,u_{\A},\xi_{\A},(a_{\A\B})_{r})$ is the induced
structure by the structure $\PP$ on $(\MM,\G)$. We denoted by $\n$
the induced Levi-Civita connection on M. We assume that
$(N_{1},...,N_{r}):=(N_{\A})$ is an orthonormal basis in the
normal space $T_{x}M^{\bot}$ at M in every point $x \in M$. In the
following, we shall identify the vector fields in M and their
images under the differential mapping, that is, if we denote the
immersion of M in N by $i$ and X is a vector field in M, we
identify X and $i_{*}X$, for all $X \in \chi(M)$.

 The Gauss and Weingarten formulae are:
\[ \widetilde{\n}_{X}Y=\n_{X}Y+\sum_{\A=1}^{r}h_{\A}(X,Y)N_{\A},\leqno(2.1)\]
and
\[ \widetilde{\n}_{X}N_{\A}=-A_{\A}X+\n_{X}^{\bot}N_{\A}, \leqno(2.2)\]
respectively, where \[ h_{\A}(X,Y)=g(A_{\A}X,Y), \leqno(2.3)\] for
every $X, Y \in \chi(M)$.

For the normal connection $\n_{X}^{\bot}N_{\A}$, we have the
decomposition
\[ \n_{X}^{\bot}N_{\A}=\sum_{\B=1}^{r}l_{\A\B}(X)N_{\B},
\leqno(2.4)\]for every $X \in \chi(M)$. Therefore, we obtain an $r
\times r$ matrix $(l_{\A\B}(X))_{r}$ of 1-forms on M. From
$\G(N_{\A},N_{\B})=\delta_{\A\B}$ we get
\[\G(\n_{X}^{\bot}N_{\A},N_{\B})+\G(N_{\A},\n_{X}^{\bot}N_{\B})=0\]
which is equivalent with
\[\G(\sum_{\gamma}l_{\A\gamma}(X)N_{\gamma},N_{\B})+
\G(N_{\A},\sum_{\gamma}l_{\B\gamma}(X)N_{\gamma})=0, \] for any $X
\in \chi(M)$. Thus, we have
\[ l_{\A\B}=-l_{\B\A}, \leqno(2.5)\]
for any $\A, \B \in \{1,...,r\}$.

For $\varepsilon =1$, the following formulae were demonstrated by
T.Adati (in \cite{Adati2}).

\thm Let M be an n-dimensional submanifold of codimension r in a
Riemannian manifold $(\MM, \G,\PP)$ (with the properties (1.1) and
(1.2)). If the structure $\PP$ is parallel with respect to the
Levi-Civita connection $\widetilde{\n}$ of $\G$, then the
$(P,g,u_{\A},\epsilon\xi_{\A},(a_{\A\B})_{r})$ induced structure
on the submanifold M has the following properties:
\[ \leqno(2.6)\begin{cases}
(i) \:(\n_{X}P)(Y)=\varepsilon\sum_{\A}h_{\A}(X,Y)\xi_{\A}+\sum_{\A}u_{\A}(Y)A_{\A}X,\\
(ii) \: (\n_{X}u_{\A})(Y)=-h_{\A}(X,PY)+\sum_{\B}u_{\B}(Y)l_{\A\B}(X)+\sum_{\B}h_{\B}(X,Y)a_{\B\A}\\
(iii) \: \n_{X}\xi_{\A}=-\varepsilon P(A_{\A}X)+\varepsilon
\sum_{\B}a_{\A\B}A_{\B}X+ \sum_{\B}l_{\A\B}(X)\xi_{\B},\\
(iv)\: X(a_{\A\B})=-\varepsilon u_{\A}(A_{\B}X)-u_{\B}(A_{\A}X)+
\sum_{\gamma}[l_{\A\gamma}(X)a_{\gamma \B}+ l_{\B
\gamma}(X)a_{\A\gamma}]
\end{cases}  \]
\normalfont \textbf{Proof:}
 From the assumption that $\widetilde{\n}\PP = 0 $ we have
\[ \widetilde{\n}_{U}(\PP V) = \widetilde{P}(\widetilde{\n}_{U}V), \leqno(2.7)\]
for any tangential vector fields U and V on $\MM$. Using the Gauss
and Weingarten formulae, we obtain from (2.4) that
\[ \widetilde{\n}_{X}(\PP Y) = \widetilde{\n}_{X}PY+
\sum_{\A}X(u_{\A}(Y))N_{\A}+\sum_{\A}u_{\A}(Y)\widetilde{\n}_{X}N_{\A}=\leqno(2.8)\]
\[=\n_{X}PY-\sum_{\A}u_{\A}(Y)A_{\A}X+\sum_{\A}[h_{\A}(X,PY)+X(u_{\A}(Y))
+\sum_{\B}u_{\B}(Y)l_{\B\A}(X)]N_{\A}\]

On the other hand, we have
\[ \PP(\widetilde{\n}_{X}Y)=\PP(\n_{X}Y)+\sum_{\A}h_{\A}(X,Y)\PP N_{\A}\] and from this we obtain

\[\PP(\widetilde{\n}_{X}Y)=P(\n_{X}Y)+\varepsilon \sum_{\A}h_{\A}(X,Y)\xi_{\A}+\leqno(2.9)\]
\[+\sum_{\A}[u_{\A}(\n_{X}Y)+\sum_{\B}h_{\B}(X,Y)a_{\B\A}]N_{\A} \]

From \cite{BYC}, we know that:
\[(\n_{X}P)(Y)=\n_{X}(PY)-P(\n_{X}Y) \leqno(2.10)\] and
\[(\n_{X}u_{\A})(Y)=X(u_{\A}(Y))-u_{\A}(\n_{X}Y) \leqno(2.11)\]

Using the relations (2.8), (2.9) in (2.7), we obtain (i) and (ii)
from (2.6), from the equality of the tangential components of M
(and the normal components of M, respectively) from the both parts
of the equality (2.7).

In the next, we apply $\widetilde{\n}_{X}$ in (1.5) (with $X \in
\chi(M)$) and using the equality (2.7), (with $Y=N_{\A}$), we get
\[ \widetilde{\n}_{X}(\PP N_{\A})=\PP(\widetilde{\n}_{X}N_{\A})\leqno(2.12)\]
From (2.12) we obtain
\[ \widetilde{\n}_{X}(\PP N_{\A})=\widetilde{\n}_{X}(\varepsilon
\xi_{\A}+\sum_{\B}a_{\A\B}N_{\B})=\leqno(2.13)\]
\[=\varepsilon \n_{X}\xi_{\A} -\sum_{\B}a_{\A\B}A_{\B}X+
\sum_{\B}[X(a_{\A\B})+\varepsilon
h_{\B}(X,\xi_{\A})+\sum_{\gamma}a_{\A\gamma}\cdot l_{\gamma
\B}(X)]N_{\B} \] and
\[\PP(\widetilde{\n}_{X}N_{\A})=\PP(-A_{\A}X+\sum_{\B}l_{\A\B}N_{\B})=\leqno(2.14)\]
\[=-P(A_{\A}X)+\varepsilon \sum_{\B}l_{\A\B}(X)\xi_{\B}-
\sum_{\B}[u_{\B}(A_{\A}X)-\sum_{\gamma}a_{\gamma\B}l_{\A\gamma}(X)]N_{\B}\]

Using the relations (2.13) and (2.14) in the equality  (2.12) and
identifying the tangential and the normal components at M,
respectively, we obtain the relations (iii) and (iv) from (2.6).
$\square $

\rem \normalfont The compatibility condition
$\widetilde{\n}\widetilde{P}=0$, where $\widetilde{\n}$ is
Levi-Civita connection with respect of the metric $\G$ implies the
integrability of the structure $\widetilde{P}$ which is equivalent
with the vanishing of the Nijenhuis torsion tensor field  of
$\widetilde{P}$:
\[N_{\widetilde{P}}(X,Y)=[\widetilde{P}X,\widetilde{P}Y]+\widetilde{P}^{2}[X,Y]-
\widetilde{P}[\widetilde{P}X,Y]-\widetilde{P}[X,\widetilde{P}Y].\leqno(2.15)\]

For this assumption, we have the next general lemma:

\lem We suppose that we have an almost product structure  $Q$ on a
manifold M and a linear connection D with the torsion T. If
$N_{Q}$ is Nijenhuis torsion tensor field  of Q, then we obtain:
\[ N_{Q}(X,Y)=(D_{QX}Q)(Y))-(D_{QY}Q)(X)+(D_{X}Q)(QY)-\leqno(2.16)\]
\[-(D_{Y}Q)(QX)-T(QX,QY)-T(X,Y)+QT(QX,Y)+QT(X,QY),\]
for any $X,Y \in \chi(M)$.\\
\normalfont \textbf{Proof:} From the definition of the torsion T
follows that :
\[[X,Y]=D_{X}Y-D_{Y}X-T(X,Y), \leqno(2.17)\]and from this we get
\[[QX,QY]=D_{QX}QY-D_{QY}QX-T(QX,QY),\leqno(2.18)\]
and
\[[QX,Y]=D_{QX}Y-D_{Y}QX-T(QX,Y),\leqno(2.19)\]
and
\[[X,QY]=D_{X}QY-D_{QY}X-T(X,QY). \leqno(2.20)\]
Using the relations $(D_{X}Q)Y=D_{X}QY-Q(D_{X}Y)$, $Q^{2}=I$ and
$D_{X}Y=D_{X}Q^{2}Y=D_{X}Q(QY)$ and replacing the relations
(2.17), (2.18), (2.19) and (2.20) in the formula of Nijenhuis
torsion tensor field of Q, we obtain :
\[ N_{Q}(X,Y)=(D_{QX}Q)(Y)-(D_{QY}Q)(X)+D_{X}Y-D_{Y}X+Q(D_{Y}QX)-\]
\[-Q(D_{X}QY)-T(QX,QY)-T(X,Y)+QT(QX,Y)+QT(X,QY)=\]
\[=(D_{QX}Q)(Y)-(D_{QY}Q)(X)+(D_{X}Q(QY)-Q(D_{X}QY))-(D_{Y}Q(QX)-\]\[-Q(D_{Y}QX))-
T(QX,QY)-T(X,Y)+QT(QX,Y)+QT(X,QY)\] so, we obtain the equality
(2.16). $\square$

\rem \normalfont From the last mentioned lemma, we remark that if
we have $T=0$ and $DQ=0$ then the structure Q is integrable. An
integrable almost product structure is also called locally product
structure.

\cor If $M$ is a totally geodesic submanifold in a locally product
manifold $(\widetilde{M},\widetilde{P},\widetilde{g})$ and the
normal connection $\n^{\bot}$ vanishes identically (that is
$l_{\A\B}=0$), then the
 $(P,g,u_{\A},\xi_{\A},(a_{\A\B})_{r})$ induced structure on M has the following properties:
\[\n P=0, \quad \n u=0, \quad \n \xi=0, \quad a=constant \]

\normalfont In the following issue, we suppose that $(\MM,\G,\PP)$
is an almost product Riemannian manifold, endowed by a linear
connection $\widetilde{\n}$ such that $\widetilde{\n}\PP =0$ and
with the torsion $\widetilde{T} \neq 0$. Let (M,g) be a
submanifold of the
 $(\MM,\G,\PP)$ Riemannian manifold, endowed with the linear metric g induced
on M by the metric $\G$ and let $\n$ be the induced connection on
M by the connection $\widetilde{\n}$ of $\MM$ .

The Gauss formula has the usual form:
\[\widetilde{\n}_{X}Y=\n_{X}Y+h(X,Y) \: (\forall) X,Y \in \chi(M)\]
where \[h(X,Y)=\sum_{\A}h^{\A}(X,Y)N_{\A}\] and $\n$ is a metric
connection on M (i.e $\n g=0$) but it is not the Levi-Civita
connection of g and his torsion has the form
\[T(X,Y)=\widetilde{T}(X,Y)-h(X,Y)+h(Y,X),\leqno(2.21)\] for any $X,Y \in
\chi(M)$. We remark that the second fundamental  form h is
bilinear in X and Y, but it is not symmetric. If the torsion
$\widetilde{T}=0$, then $T=0$ if and only if the second
fundamental form h is symmetric.

Let $c:[a,b] \longrightarrow M$, $t \longmapsto c(t)$ a smooth
curve on M. We denote by \[ \dot{c}: t \longmapsto \dot{c}(t)=
\frac{dx^{i}}{dt} \frac{\p}{\p x^{i}}_{/c(t)}\leqno(2.22)\] the
tangent vector field to c.

\defn If $\n_{\dot{c}}\dot{c}=0$ then we say that the
curve c is M-autoparallel.

\prop If c is $\MM$-autoparallel curve, then c is also,
M-autoparallel curve and $h(\dot{c},\dot{c})=0$.\\
\normalfont \textbf{Proof:} From the Gauss formula for
$X=Y=\dot{c}$ we get
\[\widetilde{\n}_{\dot{c}}\dot{c}=\n_{\dot{c}}\dot{c}+h(\dot{c},\dot{c}) \leqno(2.23)\]
But c is $\MM$-autoparallel, so
$\widetilde{\n}_{\dot{c}}\dot{c}=0$ and using the equality (2.23)
we obtain  $\n_{\dot{c}}\dot{c}=0$ and $h(\dot{c},\dot{c})=0$.
$\square$

\defn A submanifold M is said to be autoparallel in $\MM$ if any M-autoparallel curve
of submanifold M in $\MM$ is also $\MM$-autoparallel.

\normalfont We denote by
\[s_{h}(X,Y)=\frac{1}{2}(h(X,Y)+h(Y,X))\] the symmetric part and
by
\[a_{h}(X,Y)=\frac{1}{2}(h(X,Y)-h(Y,X))\] the skew-symmetric part,
respectively, of the bilinear form h. We remark that
\[h(X,Y)=s_{h}(X,Y)+a_{h}(X,Y),\leqno(2.24)\]
for any $ X,Y \in \chi(M)$.

\prop A submanifold $M \subset \MM$ is autoparallel in $\MM$ if
and only if $s_{h}=0$.\\
\normalfont \textbf{Proof:} From the Gauss formula for
$X=Y=\dot{c}$ we obtain
\[\widetilde{\n}_{\dot{c}}\dot{c}=\n_{\dot{c}}\dot{c}+s_{h}(\dot{c},\dot{c}) \leqno(2.25)\]
If $s_{h}=0$ follows that
$\widetilde{\n}_{\dot{c}}\dot{c}=\n_{\dot{c}}\dot{c}$, so any
M-autoparallel curve is also $\MM$-autoparallel.

Conversely, if any M-autoparallel curve is also
$\MM$-autoparallel, then $s_{h}(X,X)=0$ for any $X \in \chi(M)$.
Particularly, we have $s_{h}(X+Y,X+Y)=0$ and from this we obtain
$s_{h}(X,Y)=0$ for any $X,Y \in \chi(M)$, so $s_{h}=0$. $\square$

\rem \normalfont The Weingarten formula is not affected by the
non-vanishing of the torsion $\widetilde{T}$ on $\MM$ and of the
torsion T on M, thus
\[\widetilde{\n}_{X}N=-A_{N}X +\n^{\bot}_{X}\xi, \] for any $ X \in \chi(M)$ and $N \in
\Gamma(TM^{\bot})$. If $Y \in \chi(M)$ we have $\G(Y, N)=0$ so,
$X(\G(Y,N))=0$ for any $X,Y \in \chi(M)$ and this equality is
equivalent with
\[\G(\widetilde{\n}_{X}Y,N)+\G(Y,\widetilde{\n}_{X}N)=0 \leqno(2.26) \]
Using the Gauss and Weingarten formulae in (2.26), we obtain
\[g(A_{N}X,Y)=\G(h(X,Y),N),  \leqno(2.27)\]
and
\[g(A_{N}Y,X)=\G(h(Y,X),N),  \leqno(2.27)'\]
for any $X,Y \in \chi(M), N \in \Gamma(TM^{\bot})$. Thus, from
(2.27) and (2.27)' we get
\[g(A_{N}Y,X)+g(A_{N}X,Y)=2\G(s_{h}(X,Y),N), \leqno(2.28)\]
for any $X,Y \in \chi(M), N \in \Gamma(TM^{\bot})$.

If $s_{h}=0$ then we obtain $g(A_{N}Y,X)=-g(A_{N}X,Y)$.

\prop A submanifold $M \subset \MM$ is autoparallel in $\MM$ if
and only if
\[g(A_{N}Y,X)+g(A_{N}X,Y)=0, \leqno(2.29)\] for any $ X,Y \in \chi(M)$ and $N \in
\Gamma(TM^{\bot})$.

\normalfont In the following statements, we suppose that
$\widetilde{\n}_{X}\PP Y \neq \widetilde{P}(\widetilde{\n}_{X}Y)$.

We denoted by $\mathcal{P}$ the (1,2)-tensor field on $\MM$, such
that
\[\mathcal{P}(X,Y)=\widetilde{\n}_{X}\PP Y - \widetilde{P}(\widetilde{\n}_{X}Y), \leqno(2.30)\] and
\[\mathcal{P}(X,N)=\widetilde{\n}_{X}\PP N - \widetilde{P}(\widetilde{\n}_{X}N), \leqno(2.31)\]
for any $ X \in \chi(M)$ and $N \in \Gamma(TM^{\bot})$.

We denote the tangential and  normal components on M of
$\mathcal{P}(X,Y)$ by $\mathcal{P}(X,Y)^{\top}$ and
$\mathcal{P}(X,Y)^{\bot}$ , respectively, and the tangential and
normal components on M of $\mathcal{P}(X,N_{\A})$ by
$\mathcal{P}(X,N_{\A})^{\top}$ and $\mathcal{P}(X,N_{\A})^{\bot}$,
respectively.

If we omit to put the condition $\widetilde{\n}\PP=0$ in the
Theorem 2.1, then we obtain a generalization of this:

\thm Let M be an n-dimensional submanifold of codimension r in a
Riemannian manifold $(\MM,\G,\PP)$ (which satisfied the conditions
(1.1) and (1.2)). Then the structure
$(P,g,u_{\A},\varepsilon\xi_{\A},(a_{\A\B})_{r})$ induced
 on M by the structure $\PP$ has the following properties:
\[ \leqno(2.32) \normalsize\begin{cases}
(i) (\n_{X}P)(Y)=\mathcal{P}(X,Y)^{\top}+\varepsilon\sum_{\A}h_{\A}(X,Y)\xi_{\A}+\sum_{\A}u_{\A}(Y)A_{\A}X,\\
(ii)(\n_{X}u_{\A})(Y)=\G(\mathcal{P}(X,Y),N_{\A})-h_{\A}(X,PY)+\\
 \hspace{1in} +\sum_{\B}(u_{\B}(Y)l_{\A\B}(X)+h_{\B}(X,Y)a_{\B\A})\\
(iii)  \n_{X}\xi_{\A}=\mathcal{P}(X,N_{\A})^{\top}-\varepsilon
P(A_{\A}X)+\varepsilon
\sum_{\B}a_{\A\B}A_{\B}X+\sum_{\B}l_{\A\B}(X)\xi_{\B},\\
(iv)X(a_{\A\B})=\G(\mathcal{P}(X,N_{\A}),N_{\B})-\varepsilon
u_{\A}(A_{\B}X)-u_{\B}(A_{\A}X)+\\ \hspace{1in} +
\sum_{\gamma}[l_{\A\gamma}(X)a_{\gamma \B}+ l_{\B
\gamma}(X)a_{\A\gamma}]
\end{cases}
 \] for any $X,Y \in \chi(M)$.\\
\normalfont \textbf{Proof:} From (2.8) we have
\[ \widetilde{\n}_{X}(\PP Y) = \n_{X}PY-\sum_{\A}u_{\A}(Y)A_{\A}X+\]
\[+\sum_{\A}[h_{\A}(X,PY)+X(u_{\A}(Y))+\sum_{\B}u_{\B}(Y)l_{\B\A}(X)]N_{\A}\] and from (2.9) we have
\[ \PP(\widetilde{\n}_{X}Y)=P(\n_{X}Y)+\varepsilon\sum_{\A}h_{\A}(X,Y)\xi_{\A}+\sum_{\A}[u_{\A}(\n_{X}Y)
+\sum_{\B}h_{\B}(X,Y)a_{\B\A}]N_{\A}\] From the last two
equalities we obtain
\[\mathcal{P}(X,Y)=(\n_{X}P)(Y)-\sum_{\A}u_{\A}(Y)A_{\A}X-\varepsilon\sum_{\A}h_{\A}(X,Y)\xi_{\A}+ \leqno(2.33)\]
\[+\sum_{\A}[h_{\A}(X,PY)+(\n_{X}u_{\A})(Y)+\sum_{\B}u_{\B}(Y)l_{\B\A}(X)
-\sum_{\B}h_{\B}(X,Y)a_{\B\A}]N_{\A}.\]  Thus, identified the
tangent and the normal parts respectively, from (2.33) we obtain
(i) and (ii) from (2.32).

From (2.13) we have
\[ \widetilde{\n}_{X}(\PP N_{\A})=\varepsilon \n_{X}\xi_{\A}
-\sum_{\B}a_{\A\B}A_{\B}X+\sum_{\B}[X(a_{\A\B})+\varepsilon
h_{\B}(X,\xi_{\A})+\sum_{\gamma}a_{\A\gamma}\cdot l_{\gamma
\B}(X)]N_{\B}\] and from (2.14) we have
\[\PP(\widetilde{\n}_{X}N_{\A})=-P(A_{\A}X)+\varepsilon
\sum_{\B}l_{\A\B}(X)\xi_{\B}-\sum_{\B}[u_{\B}(A_{\A}X)-\sum_{\gamma}a_{\gamma\B}l_{\A\gamma}(X)]N_{\B}\]
Replacing the last two equalities in (2.31) we obtain
\[ \mathcal{P}(X,N)=\varepsilon \n_{X}\xi_{\A}+P(A_{\A}X)-\varepsilon \sum_{\B}l_{\A\B}(X)\xi_{\B}
-\sum_{\B}a_{\A\B}A_{\B}X+\leqno(2.34)\]
\[+\sum_{\B}[X(a_{\A\B})+\varepsilon h_{\B}(X,\xi_{\A})+u_{\B}(A_{\A}X)-\sum_{\gamma}
(a_{\gamma\B}l_{\A\gamma}(X)-a_{\A\gamma}\cdot l_{\gamma
\B}(X))]N_{\B}\]  Identifying the tangential and normal
components, respectively, of
$\mathcal{P}(X,N_{\A})$ from (2.34), we obtain the relations (iii) and (iv) from (2.32). $ \square $  \\


\section{\bf The normality conditions of the $(P,g,u_{\A},\xi_{\A},(a_{\A\B})_{r})$ structure }

Let M be an n-dimensional submanifold of codimension r in a
Riemannian almost product manifold $(\MM,\G,\PP)$. We suppose that
the $(P,g,u_{\A},\xi_{\A},(a_{\A\B})_{r})$ induced structure on M,
is an $(a,1)f$ Riemannian structure, where the elements
$P,g,u_{\A}$, $\xi_{\A},(a_{\A\B})_{r}$ were defined in section 1.
Let $\widetilde{\n}$ and $\n$ be the Levi-Civita connections
defined on M and $\MM$ respectively, with respect to $\G$ and g
respectively.

The Nijenhuis torsion tensor field of P has the form
\[ N_{P}(X,Y)=[PX,PY]+P^{2}[X,Y]-P[PX,Y]-P[X,PY], \leqno(3.1)\]
for any $X,Y \in \chi(M)$.

As in the case of an almost paracontact structure
(\cite{IMihai12}), one can defined the normal
$(P,g,u_{\A},\xi_{\A},(a_{\A\B})_{r})$ structure on M.

\defn If we have the equality
\[ N_{P}(X,Y)-2\sum_{\A}du_{\A}(X,Y)\xi_{\A}=0, \leqno(3.2)\] for any $X,Y\in
\chi(M)$, then the $(P,g,u_{\A},\xi_{\A},(a_{\A\B})_{r})$ induced
structure on submanifold M in a Riemannian almost product manifold
$(\MM,\G,\PP)$ is said to be normal.

\normalfont First of all, we mention a general proposition:
 \prop If $(M,g,P)$ is an almost product manifold, then the Nijenhuis
 tensor of P verifies that
\[ \quad  N_{P}(X,Y)=(\n_{PX}P)(Y)-(\n_{PY}P)(X)- \leqno(3.3)\]\[-P[(\n_{X}P)(Y)-(\n_{Y}P)(X)],  \]
for any $X,Y \in \chi(M)$, where $\n$ is the Levi-Civita connection on M .\\
\normalfont \textbf{Proof:} From the equality
\[(\n_{X}P)(Y)=\n_{X}(PY)-P(\n_{X}Y),\leqno(3.4)\] for any $X,Y
\in \chi(M)$, we obtain
\[\n_{X}(PY)=(\n_{X}P)(Y)+P(\n_{X}Y),\leqno(3.5)\]for any $X,Y \in \chi(M)$. Inverting X with Y in (3.5) we obtain
\[\n_{Y}(PX)=(\n_{Y}P)(X)+P(\n_{Y}X), \leqno(3.5)'\]for any $X,Y \in
\chi(M)$. So, from the relations (3.5) and (3.5)', it follows that
\[\n_{X}(PY)-\n_{Y}(PX)=(\n_{X}P)(Y)-(\n_{Y}P)(X)+P([X,Y]), \leqno(3.6) \]
for any $X,Y \in \chi(M)$. Then, inverting X with PX in the
equality (3.4) we obtain
\[(\n_{PX}P)(Y)=\n_{PX}(PY)-P(\n_{PX}Y),\leqno(3.7)\]for any $X,Y \in \chi(M)$. Inverting X with Y in the last relation,
 we obtain
\[(\n_{PY}P)(X)=\n_{PY}PX-P\n_{PY}X,\: (\forall) X,Y \in \chi(M) \leqno(3.7)'\]
Using (3.7) and (3.7)' in  \[ [PX,PY]=\n_{PX}PY-\n_{PY}PX\] it
follows that
\[ [PX,PY]=(\n_{PX}P)(Y)-(\n_{PY}P)(X)+P(\n_{PX}Y-\n_{PY}X), \leqno(3.8) \]
for any $ X,Y \in \chi(M)$. On the other hand, we have
\[P[PX,Y]=P\n_{PX}Y-P\n_{Y}PX,\: (\forall) X,Y \in \chi(M)
\leqno(3.9)\]and
\[P[X,PY]=P\n_{X}PY -P\n_{PY}X,\: (\forall) X,Y \in \chi(M) \leqno(3.10) \]
From $(3.8), (3.9), (3.10)$ we obtain
\[ N_{P}(X,Y)=(\n_{PX}P)(Y)-(\n_{PY}P)(X)
+P\n_{PX}Y-P\n_{PY}X+P^{2}([X,Y])-\]\[-P\n_{PX}Y+P\n_{PY}X-P[\n_{X}(PY)-\n_{Y}(PX)]\]
and using (3.6) in the last equality we have
\[N_{P}(X,Y)=(\n_{PX}P)(Y)-(\n_{PY}P)(X)+P^{2}([X,Y])-\]
\[-P[(\n_{X}P)(Y)-(\n_{Y}P)(X)]-P^{2}([X,Y])\]
and from this we obtain the equality (3.3). $\square $

From \cite{Matsumoto} we have:
\defn If $(\MM,\G,\PP)$ is an almost product Riemannian manifold
such that $\widetilde{\n} \PP =0$, then we say that $\MM$ is a
locally product Riemannian manifold.

\cor If M is a totally geodesic submanifold of a locally product
Riemannian manifold $(\MM,\G,\PP)$, with the
$(P,g,u_{\A},\xi_{\A},(a_{\A\B})_{r})$ induced structure on M,
then we have $N_{P}(X,Y)=0$, for any $X, Y \in \chi(M)$.

\thm If M is a submanifold of a locally product Riemannian
manifold $(\MM,\G,\PP)$, with
$(P,g,u_{\A},\xi_{\A},(a_{\A\B})_{r})$ induced structure and $\n$
is the Levi-Civita connection defined on M with respect to g then,
the Nijenhuis torsion tensor field of P has the form:
\[ N_{P}(X,Y)=-\sum_{\A}g((PA_{\A}-A_{\A}P)(X),Y)\xi_{\A}-\leqno(3.11)\]
\[-\sum_{\A}g(Y,\xi_{\A})(PA_{\A}-A_{\A}P)(X)+\sum_{\A}g(X,\xi_{\A})(PA_{\A}-A_{\A}P)(Y)\]
for any $X,Y \in \chi(M)$.\\
\normalfont \textbf{Proof:} From (2.6)(i) (Theorem 2.1) for
$\varepsilon = 1$, we obtain
\[(\n_{X}P)(Y)=\sum_{\A}g(A_{\A}X,Y)\xi_{\A}+\sum_{\A}g(Y,\xi_{\A})A_{\A}X \leqno(3.12)\]
and if we invert X by Y, we obtain
\[(\n_{Y}P)(X)=\sum_{\A}g(A_{\A}Y,X)\xi_{\A}+\sum_{\A}g(X,\xi_{\A})A_{\A}Y \leqno(3.12)'\]
Replacing X with PX in the equality (3.12), we obtain
\[ (\n_{PX}P)(Y)=\sum_{\A}g(A_{\A}PX,Y)\xi_{\A}+\sum_{\A}g(Y,\xi_{\A})A_{\A}(PX) \leqno(3.13)\]
If we invert X by Y then we obtain
\[ (\n_{PY}P)(X)=\sum_{\A}g(A_{\A}PY,X)\xi_{\A}+\sum_{\A}g(X,\xi_{\A})A_{\A}PY \leqno(3.13)'\]

 Replacing the relations (3.12), (3.12)', (3.13), (3.13)' in the equality  (3.3) it
 follows that
\[N_{P}(X,Y)=\sum_{\A}[g(A_{\A}PX,Y)-g(A_{\A}PY,X)]\xi_{\A}+ \leqno(3.14)\]
\[+\sum_{\A}[g(Y,\xi_{\A})A_{\A}(PX)-g(X,\xi_{\A})A_{\A}(PY)-
(\underbrace{g(A_{\A}X,Y)-g(A_{\A}Y,X))}_{=0}P\xi_{\A}]-\]
\[-\sum_{\A}[g(Y,\xi_{\A})P(A_{\A}X)-g(X,\xi_{\A})P(A_{\A}Y)]\]
But we have
\[ g(A_{\A}PY,X)=g(PY,A_{\A}X)=g(Y,PA_{\A}X)=g(PA_{\A}X,Y) \leqno(3.15)\]
and using the equality (3.15) in (3.14) we obtain (3.11).
$\square$

\cor Let M be a submanifold of codimension r in a locally product
Riemannian manifold $(\MM,\G,\PP)$, with a
$(P,g,u_{\A},\xi_{\A},(a_{\A\B})_{r})$ induced structure on M and
let $\n$ be the Levi-Civita connection defined on M with respect
to g. If (1,1) tensor field P on M commutes with the Weingarten
operators $A_{\A}$ (that is $PA_{\A}=A_{\A}P$, for any $\A \in
\{1,...,r\}$) then, the Nijenhuis torsion tensor field of P
vanishes on M (that is $N_{P}(X,Y)=0$, for any $X,Y \in \chi(M)$).

\prop Let M be a submanifold of codimension r in a locally product
Riemannian manifold $(\widetilde{M},\G,\widetilde{P})$. If
$(P,g,u_{\A},\xi_{\A},(a_{\A\B})_{r})$ is induced structure on M
and $\n$ is the Levi-Civita induced connection on M by
$\widetilde{\n}$ from $\widetilde{M}$ then, the 1-forms $u_{\A}$
verify the equality :
\[ 2du_{\A}(X,Y)=-g((PA_{\A}-A_{\A}P)(X),Y)+ \leqno(3.16)\]
\[+\sum_{\B}[l_{\A\B}(X)g(Y,\xi_{\B})-l_{\A\B}(Y)g(X,\xi_{\B})]\]
for any $X,Y \in \chi(M)$, where $l_{\A\B}$ are
the coefficients of the normal connection in the normal bundle $T^{\bot}(M)$.\\
\normalfont \textbf{Proof:} We know that
 \[2du_{\A}(X,Y)=X(u_{\A}(Y))-Y(u_{\A}(X))-u_{\A}([X,Y]),\leqno(3.17)\]
for any $ X,Y \in \chi(M)$ and
$(\n_{X}u_{\A})(Y)=X(u_{\A}(Y))-u_{\A}(\n_{X}Y)$. Thus, we obtain
\[ X(u_{\A}(Y))=(\n_{X}u_{\A})(Y)+u_{\A}(\n_{X}Y), \: (\forall) X,Y \in \chi(M) \leqno(3.18)\]
Inverting X by Y in the last relation, we obtain
\[ Y(u_{\A}(X))=(\n_{Y}u_{\A})(X)+u_{\A}(\n_{Y}X), \: (\forall) X,Y \in \chi(M) \leqno(3.18)'\]

From the relations (3.17), (3.18) and (3.18)' we have
\[ 2du_{\A}(X,Y)=(\n_{X}u_{\A})(Y)-(\n_{Y}u_{\A})(X)
+u_{\A}(\underbrace{\n_{X}Y-\n_{Y}X-[X,Y]}_{=0}),\] for any $X,Y
\in \chi(M)$ so
\[  2du_{\A}(X,Y)=(\n_{X}u_{\A})(Y)-(\n_{Y}u_{\A})(X),
\leqno(3.19)\]for any $X,Y \in \chi(M)$. From (2.6)(ii) we have
\[ (\n_{X}u_{\A})(Y)=-g(A_{\A}X,PY)+\sum_{\B}g(A_{\B}X,Y)a_{\A
\B}+\sum_{\B}g(Y,\xi_{\B})l_{\A\B}(X)\leqno(3.20)\] and inverting
X by Y in (3.20), it follows that
\[ (\n_{Y}u_{\A})(X)=-g(A_{\A}Y,PX)+\sum_{\B}g(A_{\B}Y,X)a_{\A
\B}+\sum_{\B}g(X,\xi_{\B})l_{\A\B}(Y)\leqno(3.20)'\] Replacing
(3.20) and (3.21) in the equality  (3.19) we have
\[  2du_{\A}(X,Y)=-[g(A_{\A}X,PY)-g(A_{\A}Y,PX)]+ \leqno(3.21)\]
\[+\sum_{\B}a_{\A\B}[g(A_{\B}X,Y)-g(A_{\B}Y,X)]+
\sum_{\B}[l_{\A\B}(X)g(Y,\xi_{\B})-l_{\A\B}(Y)g(X,\xi_{\B})].\]
Furthermore, we have $g(A_{\B}X,Y)=g(A_{\B}Y,X)$ and
\[g(A_{\A}X,PY)-g(A_{\A}Y,PX)=g((PA_{\A}-A_{\A}P)(X),Y)\]
and replacing the last two relations in (3.21), we obtain (3.16).
$\square$

\cor Under the assumptions of the last proposition, if the normal
connection of M vanishes identically (i.e. $l_{\A \B}=0$) then a
necessary and sufficient condition for $d u_{\A}=0$ is the
commutativity between P and $A_{\A}$ (that is $PA_{\A}=A_{\A}P$,
for any $\A \in \{1,...,r\}$).

\prop If M is a submanifold of codimension r in a locally product
Riemannian manifold $(\widetilde{M},\G,\widetilde{P})$, with the
$(P,g,u_{\A},\xi_{\A},(a_{\A\B})_{r})$ induced structure on M and
$\n$ is the Levi-Civita connection defined on M with respect to g,
then we have
\[  N_{P}(X,Y)-2\sum_{\A}du_{\A}(X,Y)\xi_{\A}=\sum_{\A}g(X,\xi_{\A})(PA_{\A}-A_{\A}P)(Y)- \leqno(3.22)\]
\[-\sum_{\A}g(Y,\xi_{\A})(PA_{\A}-A_{\A}P)(X)+\sum_{\A,\B}(g(X,\xi_{\B})l_{\A\B}(X)-g(Y,\xi_{\B})l_{\A\B}(Y))\xi_{\A},\]
for any $ X,Y \in \chi(M)$.\\
\normalfont \textbf{Proof}: The equality  (3.22) is obtained
 from (3.11) and (3.16). $\square$

\cor Let M be a submanifold of codimension r in a Riemannian
locally product manifold $(\MM,\G,\PP)$, with a
$(P,g,u_{\A},\xi_{\A},(a_{\A\B})_{r})$ induced structure on M. If
the normal connection of M is vanishes identically (i.e. $l_{\A
\B}=0$) and the (1,1) tensor field P on M commutes with the
Weingarten operators $A_{\A}$ (that is $PA_{\A}=A_{\A}P$, for any
$\A \in \{1,...,r\}$) then, M has a normal $(a,1)f$ Riemannian
structure.

\cor Let M be a submanifold of codimension r in a Riemannian
locally product  manifold $(\widetilde{M},\G,\widetilde{P})$. If
the induced structure $(P,g,u_{\A},\xi_{\A},(a_{\A\B})_{r})$ on M
is normal then we obtain
\[\sum_{\A}[g(X,\xi_{\A})(PA_{\A}-A_{\A}P)(Y)-g(Y,\xi_{\A})(PA_{\A}-A_{\A}P)(X)]+
\leqno(3.23) \]
\[+\sum_{\A,\B}[g(X,\xi_{\B})l_{\A\B}(X)-g(Y,\xi_{\B})l_{\A\B}(Y)]\xi_{\A}=0, \]
for any $X,Y \in \chi(M)$.

\cor Let M be a submanifold of codimension r in a Riemannian
locally product  manifold $(\widetilde{M},\G,\widetilde{P})$. If
the induced structure $(P,g,u_{\A},\xi_{\A},(a_{\A\B})_{r})$ on M
is normal and the normal connection $\n^{\bot}$ of M vanishes
identically (that is $l_{\A \B}=0$), then we obtain the equality
\[ \sum_{\A}g(X,\xi_{\A})(PA_{\A}-A_{\A}P)(Y)=
\sum_{\A}g(Y,\xi_{\A})(PA_{\A}-A_{\A}P)(X),\leqno(3.24) \] for any
$X,Y \in \chi(M)$.

\prop Under the assumptions of the last corollary, the equality
$(3.24)$ does not depend on the choice of a basis in the normal space $T_{x}^{\bot}(M)$, for any $x \in M$.\\
 \normalfont \textbf{Proof:} If $\{N'_{\A}\}$
is another basis in $T_{x}^{\bot}(M)$, then we have
\[N'_{\A}=\sum_{\B}s_{\A\B}N_{\B} \leqno(3.25)\]
where $(s_{\A\B})_{r}$ is an orthogonal matrix.

From the condition $\widetilde{\n}_{X}N'_{\A}=0$ we obtain
$\sum_{\B}X(s_{\A\B})N_{\B}=0$ for any $X \in M$, thus $s_{\A\B}$
are a constant functions on M. On the other hand,
\[\widetilde{\n}_{X}N'_{\A}=-A'_{\A}X \leqno(3.26)\] and
\[ \widetilde{\n}_{X}N'_{\A}=\sum_{\B}X(s_{\A\B})N_{\B}-s_{\A\B}A_{\B}X \leqno(3.27)\]
Thus, from the relations (3.25), (3.26) and (3.27) we obtain
\[A'_{\A}X=\sum_{\B}s_{\A\B}A_{\B}X \leqno(3.28)\]
Therefore, we have
\[ \PP N'_{\A}=\varepsilon \xi'_{\A}+\sum_{\B}a'_{\A \B}N_{\B}= \varepsilon \xi'_{\A}
+\sum_{\B,\C}a'_{\A\B}s_{\B\gamma}N_{\gamma} \leqno(3.29)\] and
\[ \PP N'_{\A}=\sum_{\B}s_{\A\B}\PP N_{\B}=\varepsilon\sum_{\B}s_{\A\B}\xi_{\B}+
\sum_{\B, \gamma}s_{\A\B}a_{\B\gamma}N_{\gamma} \leqno(3.30)\]So,
from (3.29) and (3.30) we obtain
\[ \xi'_{\A}=\sum_{\B}s_{\A\B}\xi_{\B} \leqno(3.31)\] and
\[\sum_{\B,\A}a'_{\A\B}s_{\B\gamma}=\sum_{\B,\gamma}s_{\A\B}a_{\B\gamma} \leqno(3.32).\]
O the basis $\{ N'_{1},...,N'_{r}\}$, the condition (3.24) becomes
\[\sum_{\A}g(X,\xi'_{\A})\cdot(PA'_{\A}-
A'_{\A}P)(Y)=\sum_{\A}g(Y,\xi'_{\A})(PA'_{\A}-A'_{\A}P)(X)
\leqno(3.33)\] From (3.28) and (3.31), we obtain
\[\sum_{\A}g(X,s_{\A\B}\xi_{\B})(Ps_{\A\gamma}A_{\gamma}-
s_{\A\gamma}A_{\gamma}P)(Y)-\leqno(3.34) \]
\[-\sum_{\A}g(Y,s_{\A\B}\xi_{\B})
(Ps_{\A\gamma}A_{\gamma}-s_{\A\gamma}A_{\gamma}P)(X)=\]
\[ =\sum_{\A}s_{\A\B}s_{\A\gamma}[g(X,\xi_{\B})\cdot(PA_{\gamma}-
A_{\gamma}P)(Y)-g(Y,\xi_{\B})(PA_{\gamma}-A_{\gamma}P)(X)]=0. \]
From the orthogonality of the matrix $(s_{\A\B})_{r}$ (that is
$\sum_{\B}s_{\A\B}s_{\gamma\B}=\delta_{\A\gamma}$) it follows that
\[\sum_{\A}[g(X,\xi_{\A})\cdot(PA_{\A}-A_{\A}P)(Y)-g(Y,\xi_{\A})(PA_{\A}-A_{\A}P)(X)]=0
\leqno(3.35)\] Therefore, the condition (3.24) does not depend on
the choice of a basis in the normal space $T_{x}^{\bot}(M)$) (for any $x \in M$). $ \square $\\

In the following we denoted by
\[ \begin{cases}
  (i) \quad B_{\A}=PA_{\A}-A_{\A}P \\
  (ii) \quad C_{\A}(X,Y)=g(B_{\A}X,Y),
\end{cases} \leqno(3.36)\]
for any  $X,Y \in \chi(M)$

\lem Let M be a submanifold in a Riemannian manifold
$(\MM,\G,\PP)$ which verifies the conditions (1.1) and (1.2). Let
$(P,g,u_{\A},\varepsilon\xi_{\A},(a_{\A\B})_{r})$ be the induced
structure on M. Then, the tensor field $C_{\A}$ on M is
skew-symmetric, thus we have
\[C_{\A}(X,Y)=-C_{\A}(Y,X), \: (\forall) X,Y \in \chi(M)\leqno(3.37)\]
\normalfont \textbf{Proof}: From (3.36)(ii) we
have:
\[C_{\A}(X,Y)=g(B_{\A}X,Y)=g(PA_{\A}X-A_{\A}PX,Y)=\]
\[=g(PA_{\A}X,Y)-g(A_{\A}PX,Y)=\]
\[=g(X,A_{\A}PY)-g(X,PA_{\A}Y)=\]
\[=-g(PA_{\A}Y-A_{\A}PY,X)=-C_{\A}(Y,X)\]
so $C_{\A}$ is skew-symmetric. $ \square $

\rem \normalfont Under the assumptions of the theorem 3.1, if we
use the notations (3.36) then the identities (3.11) and (3.16)
have the forms:
\[ N_{P}(X,Y)=\sum_{\A}[g(X,\xi_{\A})B_{\A}(Y)-
g(Y,\xi_{\A})B_{\A}(X)-C_{\A}(X,Y)\xi_{\A}]\leqno(3.38)\] and
\[2du_{\A}(X,Y)=-C_{\A}(X,Y)+
\sum_{\B}[l_{\A\B}(X)g(Y,\xi_{\B})-l_{\A\B}(Y)g(X,\xi_{\B})]\leqno(3.39)\]for
any $X,Y \in \chi(M)$.

More of them, if the induced structure on M is normal, the
equality (3.23) becomes
\[ \sum_{\A}[g(Y,\xi_{\A})B_{\A}(X)-g(X,\xi_{\A})B_{\A}(Y)]=\leqno(3.40)\]
\[=\sum_{\A,\B}[g(X,\xi_{\B})l_{\A\B}(X)-g(Y,\xi_{\B})l_{\A\B}(Y)]\xi_{\A}\]
for any $X,Y \in \chi(M)$.

\rem \normalfont Using the model for an almost paracontact
structure (\cite{IMihai12}), we can compute the components
$N^{(1)}$, $N^{(2)}$, $N^{(3)}$ and $N^{(4)}$ of the Nijenhuis
torsion tensor filed of P for the $(P,g,\xi_{\A},u_{\A},
(a_{\A\B})_{r})$ induced structure on a submanifold M in an
 almost product Riemannian manifold $(\MM,\G,\PP)$:
\[ \begin{cases}
(i) \quad N^{(1)}(X,Y)=N_{P}(X,Y)-2\sum_{\A=1}^{r}d u_{\A}(X,Y)\xi_{\A},\\
(ii) \:\: N_{\A}^{(2)}(X,Y)=(\mathcal{L}_{PX}u_{\A})Y-(\mathcal{L}_{PY}u_{\A})X,\\
(iii) \: N_{\A}^{(3)}(X)=(\mathcal{L}_{\xi_{\A}}P)X,\\
(iv) \:\: N_{\A\B}^{(4)}(X)=(\mathcal{L}_{\xi_{\A}}u_{\B})X,
\end{cases} \leqno(3.41)\]
for any $X,Y \in \chi(M)$ and $\A,\B \in \{1,...,r\}$, where
$N_{P}$ is the Nijenhuis torsion tensor field of $P$ and
$\mathcal{L}_{X}$ means the Lie derivative with respect to $X$.

\rem \normalfont Let M be a submanifold in a Riemannian almost
product manifold $(\MM,\G,\PP)$. From (3.41)(i) we remark that the
$(P,g,\xi_{\A},u_{\A}, (a_{\A\B})_{r})$ induced structure on M is
normal if and only if $N^{(1)}=0$.

\prop Let M be a submanifold in an almost product Riemannian
manifold $(\MM,\G,\PP)$, with
$(P,g,u_{\A},\xi_{\A},(a_{\A\B})_{r})$ induced structure on M. If
the normal connection $\n^{\bot}=0$ on the normal bundle
$T^{\bot}M$ vanishes identically (that is $l_{\A\B}=0$), then the
components $N^{(1)}$, $N^{(2)}$, $N^{(3)}$ and $N^{(4)}$ of the
Nijenhuis torsion tensor field of P for the structure
$(P,g,\xi_{\A},u_{\A}, (a_{\A\B})_{r})$ induced on M have the
forms:
\[\leqno(3.42)\begin{cases}
(i)\: N^{(1)}(X,Y)=\sum_{\A}g(X,\xi_{\A})(PA_{\A}-A_{\A}P)(Y)-\\
\hskip 1in -\sum_{\A}g(Y,\xi_{\A})(PA_{\A}-A_{\A}P)(X),\\
(ii)\:N_{\A}^{(2)}(X,Y)=-\sum_{\B}a_{\A\B}g((PA_{\B}-A_{\B}P)(X),Y)+\\
\hskip 0.5in +\sum_{\B}a_{\A\B}u_{\B}([X,Y])+\sum_{\B}[u_{\B}(X)u_{\A}(A_{\B}Y)-u_{\B}(Y)u_{\A}(A_{\B}X)],\\
(iii) \: N_{\A}^{(3)}(X)=\sum_{\B}a_{\A\B}(PA_{\B}-A_{\B}P)(X)-\\
\hskip 1in -P(PA_{\A}-A_{\A}P)(X)+\sum_{\B}[u_{\A}(A_{\B}X)\xi_{\B}+u_{\B}(X)A_{\B}\xi_{\A}],\\
(iv)\:
N_{\A\B}^{(4)}(X)=-u_{\A}(A_{\B}PX)-u_{\B}(PA_{\A}X)+\\
\hskip 1in
+\sum_{\gamma}[a_{\A\gamma}u_{\B}(A_{\gamma}X)+a_{\gamma
\B}u_{\A}(A_{\gamma}X)]
\end{cases} \] for any $X,Y \in \chi(M)$.\\
\normalfont \textbf{Proof:} From the equality (3.22), with the
condition $l_{\A\B}=0$, we obtain (i). Using the equality
(3.41)(ii) we have
\[N_{\A}^{(2)}(X,Y)=PX(u_{\A}(Y))-u_{\A}([PX,Y])-PY(u_{\A}(X))+u_{\A}([PY,X])=\]
\[=\underbrace{[PX(u_{\A}(Y))-u_{\A}(\n_{PX}Y)]}_{(\n_{PX}u_{\A})(Y)}-
\underbrace{[PY(u_{\A}(X))-u_{\A}(\n_{PY}X)]}_{(\n_{PY}u_{\A})(X)}+\]\[+u_{\A}(\n_{Y}PX)-u_{\A}(\n_{X}PY)=\]
\[=(\n_{PX}u_{\A})(Y)-(\n_{PY}u_{\A})(X)+\underbrace{[X(u_{\A}(PY))-u_{\A}(\n_{X}PY)]}_{(\n_{X}u_{\A})(PY)}-\]
\[-\underbrace{[Y(u_{\A}(PX))-u_{\A}(\n_{Y}PX)]}_{(\n_{Y}u_{\A})(PX)}-X(u_{\A}(PY))+Y(u_{\A}(PX))=\]
\[=\underbrace{[(\n_{PX}u_{\A})(Y)-(\n_{Y}u_{\A})(PX)]}_{2du_{\A}(PX,Y)}-
\underbrace{[(\n_{PY}u_{\A})(X)-(\n_{X}u_{\A})(PY)]}_{2du_{\A}(PY,X)}+\]
\[+X(\sum_{\B}a_{\A\B}u_{\B}(Y))-Y(\sum_{\B}a_{\A\B}u_{\B}(X))=\]
\[=2du_{\A}(PX,Y)+2du_{\A}(X,PY)+\sum_{\B}X(a_{\A\B})u_{\B}(Y)-\sum_{\B}Y(a_{\A\B})u_{\B}(X)+\]
\[+\sum_{\B}a_{\A\B}\underbrace{[X(u_{\B}(Y))-Y(u_{\B}(X))]}_{2du_{\B}(X,Y)+u_{\B}([X,Y])}\]
so we obtain
\[ N_{\A}^{(2)}(X,Y)=2du_{\A}(PX,Y)+2du_{\A}(X,PY)+
2\sum_{\B}a_{\A\B}du_{\B}(X,Y)+\leqno(3.43)\]
\[+\sum_{\B}a_{\A\B}u_{\B}([X,Y])+
\sum_{\B}u_{\B}(X)\underbrace{(u_{\A}(A_{\B}Y)+u_{\B}(A_{\A}Y))}_{-Y(a_{\A\B})}-\]
\[-\sum_{\B}u_{\B}(Y)\underbrace{(u_{\A}(A_{\B}X)+u_{\B}(A_{\A}X))}_{-X(a_{\A\B})}]\]
Using the relations (3.16) (with $l_{\A\B}=0$), (1.6)(i) and
(1.7)(iii) we obtain
\[2du_{\A}(PX,Y)+2du_{\A}(X,PY)= \]
\[=-g((PA_{\A}-A_{\A}P)(PX),Y)-g((PA_{\A}-A_{\A}P)(X),PY)=\]
\[=-g(PA_{\A}PX,Y)+g(A_{\A}P^{2}X,Y)-g(PA_{\A}X,PY)+g(A_{\A}PX,PY)\]
and from this we have
\[\leqno(3.44) 2du_{\A}(PX,Y)+2du_{\A}(X,PY)=
-\sum_{\B}u_{\B}(X)u_{\B}(A_{\A}Y)+\sum_{\B}u_{\B}(Y)u_{\B}(A_{\A}X)\]
Replacing the equality (3.44) in (3.43) we obtain
\[N_{\A}^{(2)}(X,Y)=-\sum_{\B}(u_{\B}(X)u_{\B}(A_{\A}Y)+u_{\B}(Y)u_{\B}(A_{\A}X))
-\sum_{\B}a_{\A\B}g((PA_{\B}-A_{\B}P)(X),Y)+\]
\[+\sum_{\B}(a_{\A\B}u_{\B}([X,Y])+u_{\B}(X)u_{\A}(A_{\B}Y)+u_{\B}(X)u_{\B}(A_{\A}Y)-
u_{\B}(Y)u_{\A}(A_{\B}X)-u_{\B}(Y)u_{\B}(A_{\A}X))\] and from this
we have the equality  (3.42)(ii).

Estimating  $N^{(3)}_{\A}(X)$ from the equality (3.41)(iii) we
obtain
\[N^{(3)}_{\A}(X)=\n_{\xi_{\A}}PX-\n_{PX}\xi_{\A}-P(\n_{\xi_{\A}}X)+P(\n_{X}\xi_{\A})=\]
\[= (\n_{\xi_{\A}}P)(X)-\n_{PX}\xi_{\A}+P(\n_{X}\xi_{\A}) \]
and using the relations (2.6) with the conditions $\varepsilon =1$
and $l_{\A\B}=0$ we have
\[N^{(3)}_{\A}(X)=\sum_{\B}(u_{\A}(A_{\B}X)\xi_{\B}+u_{\B}(X)A_{\B}\xi_{\A}+a_{\A\B}(PA_{\B}-A_{\B}P)(X))
-P(PA_{\A}-A_{\A}P)(X)\] and from this we obtain the equality
(iii)(3.42).

Estimating  $N^{(4)}_{\A\B}(X)$ from the equality  (3.41)(iv) we
have
\[N^{(4)}_{\A\B}(X)=\xi_{\A}(u_{\B}(X))-u_{\B}([\xi_{\A},X])=\]\[=
\underbrace{\xi_{\A}(u_{\B}(X))-u_{\B}(\n_{\xi_{\A}}X)}_{(\n_{\xi_{\A}}u_{\B})(X)}+u_{\B}(\n_{X}\xi_{\A})=
(\n_{\xi_{\A}}u_{\B})(X)+u_{\B}(\n_{X}\xi_{\A})\] and using the
relations (2.6) with the conditions $\varepsilon =1$ and
$l_{\A\B}=0$ we obtain
\[N^{(4)}_{\A\B}(X)=-h_{\B}(\xi_{\A},PX)+\sum_{\gamma}[a_{\gamma
\B}h_{\gamma}(\xi_{\A},X)+a_{\A \gamma}u_{\B}(A_{\gamma}X)]
-u_{\B}(PA_{\A}X)\] so
\[N^{(4)}_{\A\B}(X)=-g(A_{\B}PX,\xi_{\A})-g(PA_{\A}X,\xi_{\B})+\sum_{\gamma}[a_{\gamma
\B}h_{\gamma}(\xi_{\A},X)+a_{\A \gamma}u_{\B}(A_{\gamma}X)]\] and
from this we have (3.42)(iv). $\square$

\cor Under the assumptions of the last proposition, if P and the
Weingarten operators $A_{\A}$ commute (that is $PA_{\A}=A_{\A}P$,
for every $\A \in \{1,...,r\}$) then we obtain
\[\leqno(3.45)\begin{cases}
(i) \: N^{(1)}(X,Y)=0,\\
(ii)  N_{\A}^{(2)}(X,Y)=\sum_{\B}(a_{\A\B}u_{\B}([X,Y])+u_{\B}(X)u_{\A}(A_{\B}Y)-u_{\B}(Y)u_{\A}(A_{\B}X)),\\
(iii) N_{\A}^{(3)}(X)=\sum_{\B}[u_{\A}(A_{\B}X)\xi_{\B}+u_{\B}(X)A_{\B}\xi_{\A}],\\
(iv) \:
N_{\A\A}^{(4)}(X)=2\sum_{\gamma}a_{\A\gamma}u_{\A}(A_{\gamma}X)-2u_{\A}(PA_{\A}X)
\end{cases} \]


\section{\bf The commutativity of P and
$A_{\A}$ on submanifolds with
$(P,g,\xi_{\A},u_{\A},(a_{\A\B})_{r})$ normal structure}

\normalfont

Let M be an n-dimensional submanifold of codimension r in an
 almost product Riemannian manifold
$(\widetilde{M},\G,\widetilde{P})$. We suppose that M is endowed
with a $(P,g,u_{\A},\xi_{\A},(a_{\A\B})_{r})$ induced structure on
M by the $\PP$ and the normal connection $\n^{\bot}$ on the normal
bundle $T^{\bot}(M)$ vanishes identically.

First of all, we search for necessary conditions for the linearly
independent of the tangent vector fields $\xi_{1},...,\xi_{r}$
(with $r \geq 2$). In this situation, we will show that the
condition of the normality of induced structure on M is equivalent
with the commutativity between the tensor field P and the
Weingarten operators $A_{\A}$ (for $\A \in \{1,...,r\}$). We
denoted by
\[\mathcal{A}:=(a_{\A\B})_{r} \leqno(4.1)\] the matrix from the $(P,g,u_{\A},\xi_{\A},(a_{\A\B})_{r})$
induced structure on M.

\prop Let M be a submanifold of codimension r (with $r \geq 2$) in
a locally product Riemannian manifold
$(\widetilde{M},\G,\widetilde{P})$, with
$(P,g,u_{\A},\xi_{\A},(a_{\A\B})_{r})$ induced structure on M by
$\PP$. If the normal connection $\n^{\bot}$ vanishes identically
on the normal bundle $T^{\bot}(M)$ (i.e.$l_{\A\B}=0$) then, the
tangent vector fields $\{\xi_{1},...,\xi_{r}\}$ are linearly
independent if and only if the determinant of the matrix
($I_{r}-\mathcal{A}^{2}$) does not vanish in any point $x \in M$,
(where $I_{r}$ is the $r \times r$ identity matrix).\\
\normalfont \textbf{Proof:} Let $k_{1},...,k_{r}$ the real number
with the properties that
\[k_{1}\xi_{1}+...+k_{r}\xi_{r}=0 \leqno(4.2)\]
in any point $x \in M$. From the equality (1.6)(iv), for
$\varepsilon =1$, we obtain
\[g(\xi_{\A},\xi_{\B})=\delta_{\A\B}-\sum_{\gamma}a_{\A\gamma}a_{\gamma\B}
 \leqno(4.3)\]
Multiplying the equality  (4.2) by $\xi_{\A}$ (for any $\A \in
\{1,...,r\}$) and using the equality  (4.3) we obtain :
\[ \leqno(4.4) \begin{cases}
k_{1}(1-\sum_{\gamma}a_{1\gamma}a_{\gamma
1})+k_{2}(-\sum_{\gamma}a_{1\gamma}a_{\gamma
2})+...+k_{r}(-\sum_{\gamma}a_{1\gamma}a_{\gamma r})=0\\

k_{1}(-\sum_{\gamma}a_{1\gamma}a_{\gamma
2})+k_{2}(1-\sum_{\gamma}a_{2\gamma}a_{\gamma
2})+...+k_{r}(-\sum_{\gamma}a_{2\gamma}a_{\gamma r})=0\\

...............................................................................\\

k_{1}(-\sum_{\gamma}a_{1\gamma}a_{\gamma
r})+k_{2}(-\sum_{\gamma}a_{2\gamma}a_{\gamma
r})+...+k_{r}(1-\sum_{\gamma}a_{r\gamma}a_{\gamma r})=0
\end{cases} \]
This linear system of equations has the unique solution
$k_{1}=...=k_{r}=0$ if and only if it does not have a vanishing
determinant. Furthermore, the determinant of the linear system of
equations (4.4) is the determinant of the following matrix:
\[I_{r}-\begin{pmatrix}
  a_{11} & a_{12} & a_{13} & ... & a_{1r} \\
  a_{21} & a_{22} & a_{23} & ... & a_{2r} \\
  ...   & ...& ...&...&... \\
  a_{r1} & a_{r2} & a_{r3} & ... & a_{rr}
\end{pmatrix} \cdot \begin{pmatrix}
  a_{11} & a_{12} & ... & a_{1r} \\
  a_{21} & a_{22} & ... & a_{2r} \\
  a_{31} & a_{32} & ... & a_{3r} \\
  ...& ... & ... & ... \\
  a_{r1} & a_{r2} &... & a_{rr}
\end{pmatrix} \]
which is the determinant of the matrix $I_{r}-\mathcal{A}^{2}$.
$\square$

\thm Let M be an n-dimensional submanifold of codimension r in a
locally product Riemannian manifold
$(\widetilde{M},\G,\widetilde{P})$, with
$(P,g,u_{\A},\xi_{\A},(a_{\A\B})_{r})$ induced structure normal on
M. If the normal connection $\n^{\bot}$ vanishes identically
(i.e.$l_{\A\B}=0$) and the determinant of matrix
($I_{r}-\mathcal{A}^{2}$) does not vanish in any point $x \in M$,
then the (1,1) tensor field P commutes with the Weingarten
operators $A_{\A}$, thus
\[PA_{\A}=A_{\A}P, \leqno(4.5)\] for any
$X \in \chi(M)$ and $\A \in \{1,...,r\}$ \\
\normalfont \textbf{Proof:} If the normal connection vanishes
identically (thus $l_{\A\B}=0$), the equality  (3.40) can be
written in the form
\[ \sum_{\A}g(X,\xi_{\A})B_{\A}(Y)=
\sum_{\A}g(Y,\xi_{\A})B_{\A}(X),\leqno(4.6) \] for any $X,Y
\in\chi(M)$ and $\A \in \{1,...,r\}$.

Multiplying with $Z \in \chi(M)$ the equality  (4.6)
 and using (3.36)(ii) we obtain:
\[ \sum_{\A}g(X,\xi_{\A})C_{\A}(Y,Z)=
\sum_{\A}g(Y,\xi_{\A})C_{\A}(X,Z),\leqno(4.7) \] for any $X,Y \in
\chi(M)$. Inverting Y by Z in the equality (4.7), we have
\[ \sum_{\A}g(X,\xi_{\A})C_{\A}(Z,Y)=
\sum_{\A}g(Z,\xi_{\A})C_{\A}(X,Y),\leqno(4.7)' \] Summating the
relations (4.7) and (4.7)', we obtain
\[\sum_{\A}g(Y,\xi_{\A})C_{\A}(X,Z)+\sum_{\A}g(Z,\xi_{\A})C_{\A}(X,Y)=0,\leqno(4.8) \]
because $C_{\A}$ is skew symmetric (i.e.
$C_{\A}(Y,Z)+C_{\A}(Z,Y)=0$). The equality  (4.8) is equivalent
with
\[\sum_{\A}g(g(Y,\xi_{\A})B_{\A}(X)+g(B_{\A}(X),Y)\xi_{\A},Z)=0, \:
 (\forall)X,Y,Z \in \chi(M) \]
so
\[\sum_{\A}g(Y,\xi_{\A})B_{\A}(X)+\sum_{\A}g(B_{\A}(X),Y)\xi_{\A}=0,
\: (\forall)X,Y \in \chi(M) \leqno(4.9) \] Using (3.36)(ii), the
equality (4.9) can be written in the form
\[\sum_{\A}g(Y,\xi_{\A})B_{\A}(X)+\sum_{\A}C_{\A}(X,Y)\xi_{\A}=0,
\: (\forall)X,Y \in \chi(M) \leqno(4.10) \] Inverting X by Y
 in the equality (4.10) we obtain
\[\sum_{\A}g(X,\xi_{\A})B_{\A}(Y)+\sum_{\A}C_{\A}(Y,X)\xi_{\A}=0,
\: (\forall)X,Y \in \chi(M) \leqno(4.10)' \] Summating  the
equalities (4.10) and (4.10)' we obtain
\[ \sum_{\A}g(Y,\xi_{\A})B_{\A}(X)+\sum_{\A}g(X,\xi_{\A})B_{\A}(Y)=0,
\: (\forall)X,Y \in \chi(M) \leqno(4.11)\]because
$C_{\A}(Y,Z)+C_{\A}(Z,Y)=0$. Therefore, from the relations (4.6)
and (4.11) it follows that :
\[ \sum_{\A}g(Y,\xi_{\A})B_{\A}(X)=0, \: (\forall)X,Y \in \chi(M). \leqno(4.12)\]
From $det(\mathcal{A}^{2}-I_{r})\neq 0$ we obtain that the
$\xi_{1},...,\xi_{r}$ tangential vector fields are linearly
independent in any point $x \in M$, for $r \geq 2$. Because we
have r linearly independent tangent vector fields on M, then we
have $r\leq n$ and from this it follows that there exist a tangent
vector field $Y \in \chi(M)$ which is orthogonal on the space
spanned by $\{\xi_{1},...,\xi_{r}\}-\{\xi_{\A}\}$ and
$g(Y,\xi_{\A})\neq 0$. Thus, from the equality (4.12) we have
$B_{\A}(X)=0$, for any $X \in \chi(M)$ and $\A \in \{1,...,r\}$
($r>1$) and from this we have (4.5).

For $r=1$, from the equality (4.12) we have $g(Y,\xi)B(X)=0$
(where $B=PA-AP$). For $Y=\xi$ we obtain  $g(\xi,\xi)B(X)=0$. But
$g(\xi,\xi)=1-a^{2}$ and from $1-a^{2} \neq 0$ we have $B(X)=0$,
so $PA=AP$. $\square$

From the Corollary 3.5 and Theorem 4.1 we obtain:

\thm Let M be an n-dimensional submanifold of codimension r in a
locally product Riemannian manifold
$(\widetilde{M},\G,\widetilde{P})$, with
 the normal connection $\n^{\bot}$ vanishes identically
(i.e.$l_{\A\B}=0$). If the determinant of the matrix
($I_{r}-\mathcal{A}^{2}$) does not vanish in any point $x \in M$,
then the $(P,g,u_{\A},\xi_{\A},(a_{\A\B})_{r})$ induced structure
on M is normal if and only if the (1,1) tensor field P commutes
with the Weingarten operators $A_{\A}$ (for any $\A \in
\{1,...,r\}$).

\rem \normalfont Under the assumptions of the last theorem, if the
submanifold M in $\MM$ is totally umbilical (or totally geodesic),
then the commutativity between the (1,1) tensor field P and the
Weingarten operators $A_{\A}$ (for any $\A \in \{1,...,r\}$) has
done.


\section{\bf On the composition of the immersions on manifolds
with $(P,g,\xi_{\A},u_{\A}, (a_{\A\B})_{r})$-structure}

Let $(\MM,\G,\PP)$ be an  almost product Riemannian manifold and
let $(\M,\g)$ be a submanifold of codimension 1, isometric
immersed in $\MM$ (with the induced metric $\g$ on $\M$ by $\G$
and $N_{2}$ an unit vector field, normal on $\M$ in $\MM$). On the
other hand, we suppose that $(M,g)$ is an isometric immersed
submanifold of codimension 1 in $\M$ and let $N_{1}$ be an unit
vector field, normal on $M$ in $\M$. Therefore, $(M,g)$ is an
isometric immersed Riemannian submanifold in $(\MM,\G)$. We may
assume that M is imbedded in $\M$ and $\M$ is imbedded in $\MM$.

From the decomposed of the vector fields
$\widetilde{P}\overline{X}$ ( $\overline{X} \in \chi(\M)$) and
$\widetilde{P}N_{2}$ respectively, in tangential and normal
components at $\M$ in $\MM$, we obtain
\[\widetilde{P}\overline{X}=\overline{P}\:\overline{X}+u_{2}(\overline{X})N_{2}
, \leqno(5.1)\] for any $\overline{X} \in \chi(\overline{M})$ and
\[\widetilde{P}N_{2}=\xi_{2}+a_{22}N_{2}, \leqno(5.2)\]
where $\overline{P}$ is an (1.1) tensor field on $\M$, $u_{2}$ is
an 1-form on $\M$, $\xi_{2}$ is a tangent vector field on $\M$ and
$a_{22}$ is a real function on $\M$.\\

\normalfont For r=1, the relations (1.6) and (1.7) for the
submanifold $\M$ in $\MM$ are written in the next proposition:

\prop The almost product Riemannian structure $(\PP,\G)$ on a
manifold $\MM$ induces, on any submanifold $\M$ of codimension 1
in $\MM$, a $(\overline{P},\g,u_{2}, \xi_{2}, a_{22})$ Riemannian
structure (where $\overline{P}$ is an (1,1) tensor field on $\M$,
$u_{2}$ is an 1-form on $\M$, $\xi_{2}$ is a tangent vector field
on $\M$ and $a_{22}$ is a real function on $\M$) with the
following properties:
\[\begin{cases}
(i)\quad
\overline{P}^{2}\overline{X}=\overline{X}-u_{2}(\overline{X})\xi_{2},
\: (\forall) \overline{X} \in \chi(\M),\\
(ii) \quad
u_{2}(\overline{P}\:\overline{X})=-a_{22}u_{2}(\overline{X}),
\:(\forall) \overline{X} \in \chi(\M),\\
(iii)\:u_{2}(\xi_{2})=1-a_{22}^{2},\\
(iv)\: \overline{P}\xi_{2}=-a_{22}\xi_{2},\\
\end{cases} \leqno(5.3)\] and
\[\begin{cases}
(i) \quad u_{2}(\overline{X})=\g(\overline{X},\xi_{2}),
\: (\forall) \overline{X} \in \chi(\M),\\
(ii) \quad
\g(\overline{P}\:\overline{X},\overline{Y})=\g(\overline{X},\overline{P}\:\overline{Y}),
\: (\forall) \overline{X},\overline{Y} \in \chi(\M),\\
(iii)\:\g(\overline{P}\:\overline{X},\overline{P}\:\overline{Y})=
\g(\overline{X},\overline{Y}), \: (\forall)
\overline{X},\overline{Y} \in \chi(\M).\quad \square
\end{cases} \leqno(5.4) \]

\normalfont  From the decomposed of the vector fields
$\overline{P}X$ ($X \in \chi(M)$) and $\overline{P}N_{1}$
respectively, in the tangential and normal components on $M$ in
$\M$, we obtain
\[\overline{P}X=PX+u_{1}(X)N_{1}
, \quad (\forall) X \in \chi(M) \leqno(5.5)\] and
\[\overline{P}N_{1}=\xi_{1}+a_{11}N_{1}, \leqno(5.6)\]
where $P$ is an (1,1) tensor field on $M$, $u_{1}$ is an 1-form on
$M$, $\xi_{1}$ is a tangential vector field on $M$ and $a_{11}$ is
a real function on $M$.

On the other hand, the vector field $\xi_{2} \in \chi(\M)$ can be
decomposed in the tangential and normal components on M in $\M$:
\[\xi_{2}=\xi_{2}^{\top}+\xi_{2}^{\bot},\leqno(5.7)\]
and we remark that $\xi_{2}^{\bot}$ and $N_{1}$ are collinear.

 \prop  The vector fields $\widetilde{P}X$ ( $X \in
\chi(M)$), $\widetilde{P}N_{1}$ and $\widetilde{P}N_{2}$ have the
next decomposes in the tangential and normal parts on M in $\MM$:
\[\begin{cases}
(i)\quad \widetilde{P}X=PX+u_{1}(X)N_{1}+u_{2}(X)N_{2},\:(\forall) X \in \chi(M)\\
(ii) \quad \widetilde{P}N_{1}=\xi_{1}+a_{11}N_{1}+a_{12}N_{2},\\
(iii)\:\widetilde{P}N_{2}=\xi_{2}^{\top}+a_{12}N_{1}+a_{22}N_{2}
\end{cases}\leqno(5.8)
\]
where $P$ is an (1,1) tensor field on $M$, $u_{1}$ is an 1-form on
$M$, $\xi_{1}$ is a tangent vector field on $M$, $(a_{\A\B})$
(with $\A, \B \in \{1,2\}$) is an $r \times r$ matrix where its
entries $a_{11}$, $a_{22}$ and
$a_{12}=a_{21}=\G(\xi_{2}^{\bot},N_{1})$ are real functions on $M$
and $u_{2}, \xi_{2}^{\top}, a_{22}$ were defined in the last
proposition.\\
\normalfont \textbf{Proof:} From the relations (5.1)
and (5.5) we obtain
\[\PP X= \overline{P} X + u_{2}(X)N_{2}=PX+u_{1}(X)N_{1}+u_{2}(X)N_{2}\]
so, (i) has done.

 From the relations (5.1) and (5.6) we obtain
 \[\PP N_{1}= \overline{P} N_{1}+u_{2}(N_{1})N_{1}=\xi_{1}+a_{11}N_{1}+u_{2}(N_{1})N_{2}. \leqno(5.9)\]
We denote by  $u_{2}(N_{1}):=a_{12}$. Thus,
\[a_{12}=\G(\xi_{2},N_{1})=\G(\xi_{2}^{\bot},N_{1}) \leqno(5.10) \]
and from this it follows that
\[\xi_{2}^{\bot}=\G(\xi_{2},N_{1})N_{1}=u_{2}(N_{1})N_{1}=a_{12}N_{1}. \leqno(5.11)\]
So, from the equality (5.9) we have (ii) from (5.8).

From the relations (5.2) and (5.7) we have
\[ \PP N_{2}= \xi_{2}^{\top}+\xi_{2}^{\bot}+a_{22}N_{2}\]
and using (5.11) we obtain  (iii) from (5.8) (where
$a_{21}=a_{12}$). $\square$

 \thm The structure $(\overline{P},\g,u_{2}, \xi_{2}, a_{22})$
 (induced on a submanifold $(\M,\g)$ of codimension 1 in an
 almost product Riemannian manifold $(\MM,\G,\PP)$) also induces,
on a submanifold $(M,g)$ of codimension 1 in $\M$, a Riemannian
structure $(P,g,u_{1},u_{2},\xi_{1},\xi_{2}^{\top},(a_{\A\B}))$
(where $P, u_{1}, u_{2}, \xi_{1}, \xi_{2}^{\top}, (a_{\A\B})$ were
defined in the last two propositions) which has the following
properties:

\[ \leqno(5.12) \begin{cases}
(i) \quad
P^{2}X=X-u_{1}(X)\xi_{1}-u_{2}(X)\xi_{2}^{\top},\:(\forall) X \in
\chi(M)\\
(ii)\quad u_{1}(PX)=-a_{11}u_{1}(X)-a_{12}u_{2}(X),\:(\forall) X
\in \chi(M)\\
(iii)\:u_{2}(PX)=-a_{21}u_{1}(X)-a_{22}u_{2}(X),\:(\forall) X \in \chi(M)\\
(iv)\quad u_{1}(\xi_{1})=1-a_{11}^{2}-a_{12}^{2},\\
(v)\quad u_{2}(\xi_{1})=-a_{11}a_{12}-a_{12}a_{22},\\
(vi)\quad u_{1}(\xi_{2}^{\top})=-a_{11}a_{12}-a_{12}a_{22},\\
(vii)\quad u_{2}(\xi_{2}^{\top})=1-a_{12}^{2}-a_{22}^{2},\\
(viii)\:P(\xi_{1})=-a_{11}\xi_{1}-a_{12}\xi_{2}^{\top},\\
(ix)\:P(\xi_{2}^{\top})=-a_{12}\xi_{1}-a_{22}\xi_{2}^{\top},
\end{cases} \]
and the properties which depends on the metric g are:
\[ \leqno(5.13)\begin{cases}
(i)\quad u_{1}(X)=g(X,\xi_{1}),\\
(ii) \quad u_{2}(X)=g(X,\xi_{2}^{\top}),\\
(iii)\:g(PX,Y)=g(X,PY),\\
(iv)\quad g(PX,PY)=g(X,Y)-u_{1}(X)u_{1}(Y)-u_{2}(X)u_{2}(Y),
\end{cases} \]
for any $X,Y \in \chi(M).$\\
 \normalfont \textbf{Proof:} From $\PP(\PP X)=X$ and (5.8) it
follows that
\[\PP(PX+u_{1}(X)N_{1}+u_{2}(X)N_{2})=X\]
thus we have
\[P^{2}X+u_{1}(PX)N_{1}+u_{2}(PX)N_{2}+u_{1}(X)(\xi_{1}+a_{11}N_{1}+a_{12}N_{2})+\]
\[+u_{2}(X)(\xi_{2}^{\top}+a_{12}N_{1}+a_{22}N_{2})=X\]
Identifying the tangential components on M from the last equality,
we obtain (i) from (5.12). Then, multiplying the last equality by
$N_{1}$ and $N_{2}$ respectively, and using the equality (5.11) we
obtain the relations (ii) and (iii) from (5.12).

On the other hand, from $\PP(\PP N_{1})=N_{1}$ we obtain
\[N_{1}=\PP (\PP N_{1})=\PP (\xi_{1}+a_{11}N_{1}+a_{12}N_{2})\]
and using the relations (5.8) it follows that
\[N_{1}=P\xi_{1}+u_{1}(\xi_{1})N_{1}+u_{2}(\xi_{1})N_{2}+\]
\[+a_{11}(\xi_{1}+a_{11}N_{1}+a_{12}N_{2})+a_{12}(\xi_{2}^{\top}+a_{21}N_{1}+a_{22}N_{2})\]
Identifying the tangential components on M from the last equality,
and multiplying this relation by $N_{1}$ and $N_{2}$ respectively,
we obtain, from the equality  (5.11), the relations (iv), (v) and
(viii) respectively, from (5.12).

From $\PP(\PP N_{2})=N_{2}$ we obtain
\[N_{2}=\PP (\PP N_{2})=\PP (\xi_{2}^{\top}+a_{12}N_{1}+a_{22}N_{2})\]
and using the relations (5.8) it follows that
\[N_{2}=P(\xi_{2}^{\top})+u_{1}(\xi_{2}^{\top})N_{1}+u_{2}(\xi_{2}^{\top})N_{2}+\]
\[+a_{12}(\xi_{1}+a_{11}N_{1}+a_{12}N_{2})+a_{22}(\xi_{2}^{\top}+a_{21}N_{1}+a_{22}N_{2})\]
Identifying the tangential components on M from the last equality
and multiplying this relation by $N_{1}$ and  $N_{2}$
respectively, we obtain, from (5.11) the relations (vi), (vii) and
(ix) from (5.12).

From \[g(PX,Y)=\G(\PP X -u_{1}N_{1}-u_{2}N_{2},Y)=\G(\PP X,Y)=\]
\[=\G(X, \PP Y)=\G(X,PY+u_{1}(Y)+u_{2}(Y)N_{2})=g(X,PY)\]
we obtain the equality  (iii) from (5.13).

From \[\G(\PP X,N_{1})=\G (X, \PP N_{1})\]and using the relations
(5.8), we have
\[\G(PX+u_{1}(X)N_{1}+u_{2}(X)N_{2},N_{1})=\G(X,\xi_{1}+a_{11}N_{1}+a_{12}N_{2})\]
Thus, $u_{1}(X)=\G(X,\xi_{1})=g(X,\xi_{1})$ and from this we
obtain the equality (i) from (5.13).

From \[\G(\PP X,N_{2})=\G (X, \PP N_{2})\] and using the relations
(5.8) we have
\[\G(PX+u_{1}(X)N_{1}+u_{2}(X)N_{2},N_{2})=\G(X,\xi_{2}^{\top}+a_{12}N_{1}+a_{22}N_{2})\]
Thus, $u_{2}(X)=\G(X,\xi_{2}^{\top})=g(X,\xi_{2}^{\top})$ and from
this we have the equality  (ii) from (5.13).

From $g(PX,Y)=g(X,PY)$, replacing Y with PY and using the equality
(i) from  (5.12) we obtain
\[g(PX,PY)=g(X,P^{2}Y)=g(X,Y-u_{1}(Y)\xi_{1}-u_{2}(Y)\xi_{2}^{\top}).\]
Thus, from the relations (i) and (ii) from (5.13), it follows that
 \[g(PX,PY)=g(X,Y)-u_{1}(X)u_{1}(Y)-u_{2}(X)u_{2}(Y)\]
and from this, we obtain the equality  (iv) from (5.13). $\square$

\defn We say that a Riemannian immersion is an isometric immersion between two
Riemannian manifolds.

\cor Let M be a submanifold of codimension 1, isometric immersed
in $\M$, which is also of codimension 1 and isometric immersed in
an almost product Riemannian manifold $(\MM,\G,\PP)$, such that we
have the Riemannian immersions:
 \[(M,g)\hookrightarrow (\M,\g) \hookrightarrow(\MM,\G)\]
Then, the induced structure on M by the structure $(\PP,\G)$ from
$\MM$ is a $(P,g,u_{1},u_{2},\xi_{1},\xi_{2}^{\top}, (a_{\A\B}))$-
Riemannian structure  (where $P$ , $u_{1}$, $u_{2}$, $\xi_{1}$,
$\xi_{2}^{\top}$, $(a_{\A\B})$ were defined in the last theorem).
 This structure is determined by the structure $(\overline{P},\g,u_{2},\xi_{2},a_{22})$
(induced on $\M$ by the structure from $\MM$) and by the structure
$(P,g,u_{1},\xi_{1},a_{11})$ (induced on a M by the structure from
$\M$).

\cor Let $M:=M_{1}$ be a submanifold of codimension r (with $r
\geq 2 $) in an almost product Riemannian manifold $(\MM,\G,\PP)$.
We make the following notations: $\MM :=M_{r+1}$, $\G :=g^{r+1}$,
$\PP := P_{r+1}$, such that we have the sequence of Riemannian
immersions
\[(M_{1},g^{1})\hookrightarrow (M_{2},g^{2})
\hookrightarrow ...\hookrightarrow (M_{r},g^{r})\hookrightarrow
(\MM,\G,\PP)\]  where $g^{i}$ is an induced metric on $M^{i}$ by
the metric $g^{i+1}$ from $M_{i+1}$, ($i \in \{1,...,r\}$) and
each one of $(M_{i},g^{i})$ is a submanifold of codimension 1,
isometric immersed in the manifold $(M_{i+1},g^{i+1})$ ($i \in
\{1,...,r\}$). Then, the structure
$(P_{1},g^{1},\xi^{1}_{\A_{1}},u^{1}_{\A_{1}},(a^{1}_{\A_{1}\B_{1}}))$,
which is successive induced by the structures
$(P_{i},g^{i},\xi^{i}_{\A_{i}},u_{\A_{i}},(a^{i}_{\A_{i}\B_{i}}))$
on the manifolds $M_{i}$ ($i \in \{2,...,r\}$, $\A_{i}\in
\{i,...,r\}$) is the same as the induced structure on
$(M_{1},g^{1})$ by the almost product structure $\PP$ on $\MM$,
where the vectorial fields $\xi^{1}_{\A_{i}}$ on $M_{1}$ are the
tangential components at M of the tangent vectorial fields
$\xi^{i}_{\A_{i}}$ from $M_{i}$, the 1-forms $u^{1}_{\A_{i}}$ are
the restrictions on M of the 1-forms $u^{i}_{\A_{i}}$ from $M_{i}$
(for any $i \in \{2,...,r\}$ and $\A_{i}\in \{i,...,r\}$), and the
entries of the $r \times r$ matrix $(a^{1}_{\A_{1}\B_{1}})$ are
defined by
\[\begin{cases}
\quad a^{1}_{\A_{1},\A_{1}}=g^{1}(P_{1}(N_{\A_{1}}),N_{\A_{1}}), \\
\quad
a^{1}_{\A_{1}\B_{1}}=a^{1}_{\B_{1}\A_{1}}=g^{1}(\xi^{1}_{\A_{1}},N_{\B_{1}}),
\quad \A_{1} > \B_{1}
\end{cases} \leqno(5.14)\]
for any  $\A_{1},\B_{1} \in \{1,...r\}$.

\section{\bf $(P,g,\xi_{\A},u_{\A},(a_{\A \B})_{r})$ induced structure on submanifolds of codimension 1 or 2
in almost product Riemannian manifolds}

\normalfont The hypersurfaces imersed in an almost product
Riemannian manifolds were studied by T. Adati (\cite{Adati1},
\cite{Adati2}), T. Miyazawa (\cite{Miyazawa}) and M. Okumura
(\cite{Okumura}).

In the following, we assume that M is a submanifold of codimension
1 in an  almost product Riemannian manifold $(\MM,\G,\PP)$. We
denoted by $A_{N}:=A$ (where N is an unit normal vector field at
submanifold M in $\MM$), $u_{1}:=u$, $\xi_{1}:=\xi$ and
$a_{\A\B}:=a$.

 From (1.4), (1.5) we obtain
\[\PP X= PX+u(X)N,  \leqno(6.1)\]
for any $X \in \chi(M)$, and
\[\PP N = \xi+ a N \leqno(6.2)\]  where P is a tensor field on M, u is an 1-form
on M, $a$ is a real function on M and $\xi$ is a tangent vector
field on M. The induced structure on M is a $(P,g,u,\xi,a)$
Riemannian structure. From (1.6) we have:
\[  \begin{cases}
(i) \quad P^{2}X=X-u(X)\xi,\\
(ii) \:\: u(PX)=-au(X), \\
(iii) \: u(\xi)=1-a^{2},\\
(iv) \: P(\xi)=-a\xi
\end{cases} \leqno(6.3)\]
for any $X \in \chi(M)$.

From the relations (1.7) we obtain
\[ \begin{cases}
(i) \quad u(X)=g(X,\xi),\\
(ii) \:\: g(PX,PY)=g(X,Y)-u(X)u(Y)
\end{cases}\leqno(6.4)\]
for any $X,Y \in \chi(M)$.

\rem \normalfont For a=0, from the relations (1.3), we remark that
$(P,u,\xi)$ is an almost paracontact structure on M and
$(P,g,u,\xi)$ is a Riemannian almost paracontact structure.

The Gauss and Weingarten formulae on the hypersurface M, isometric
immersed in $\MM$, have the forms:
\[ \widetilde{\n}_{X}Y=\n_{X}Y + g(AX,Y)N \leqno(6.5)\]
and
\[\widetilde{\n}_{X}N = -AX \leqno(6.6)\] respectiverly,
for any $X,Y \in \chi(M)$.

If M is a  locally product Riemannian manifold, the formulae (2.6)
become:
\[ \begin{cases}
(i)\quad (\n_{X}P)(Y)=u(Y)AX+g(AX,Y)\xi,\\
(ii) \:\: (\n_{X}u)(Y)=-g(AX,PY) + a \cdot g(AX,Y),\\
(iii) \: \n_{X}\xi=-P(AX)+ a \cdot AX,\\
 (iv) \:\:  X(a)= -2u(AX)=-2g(AX,\xi)=-2g(X,A\xi)
\end{cases}\leqno(6.7)
\]
for any $X,Y \in \chi(M)$.

\rem \normalfont From the relations (6.3)(iii) and (6.4)(i), we
remark that $a \in (-1,1)$ for $\xi \neq 0$.

\rem \normalfont If $a=0$, then the relations (6.3) have the
forms:
\[  \begin{cases}
(i) \quad P^{2}X=X-u(X)\xi,\\
(ii) \:\: u(PX)=0, \\
(iii) \: u(\xi)=1,\\
(iv) \: P(\xi)=0,
\end{cases}\leqno(6.8)\]
for any $X \in \chi(M)$.

\rem \normalfont Let M be a hypersurface in an  almost product
Riemannian manifold $(\MM,\G,\PP)$. From the relations (1.3)(iii)
and (1.4)(i) we remark that $u(\xi)=0$ so, $g(\xi,\xi)=0$, for
$a^{2}=1$ and from this we have $\xi=0$. Therefore, from the
equality  (1.2) we obtain $\PP N = a \cdot N$ and $P^{2}X=X$.
Thus, the induced structure P on M becomes an almost product
structure and an eigenvalue of $\PP$ is a.

\rem \normalfont Let M be a hypersurface in an almost product
Riemannian manifold $(\MM,\G,\PP)$. From (6.3)(iv) we remark that
$\xi$ is an eigenvector of a tensor field P with the eigenvalue
$-a$.

In the following we denote by $\xi^{\bot} = \{X \in \chi(M) / X
\bot \xi\}$.

\rem \normalfont From $u(PX)=g(PX,\xi)$ and
$u(PX)=-au(X)=-ag(X,\xi)$ we obtain $g(PX+aX,\xi)=0$ and $(PX+aX)
\bot \xi, (\forall)X \in \chi(M)$.  We remark that there is a
vector $V \in \xi^{\bot}$ such that $PX=-aX+V$.

\prop Let M be a hypersurface in an almost product Riemannian
manifold $(\MM,\G,\PP)$, with the structure $(P,g,\xi,u,a)$
induced on M by the structure $\PP$ from $\MM$. We suppose that
$\xi \neq 0$. A necessary and sufficient condition for M to be
totally geodesic in $\MM$ is that $\n_{X}P=0$, for any $X \in \chi(M)$. \\
\normalfont\textbf{Proof}:
If M is totally geodesic, then $\n_{X}P=0$ ($\forall X \in \chi(M)$) from the equality (6.7)(i).\\
Conversely, we suppose that $\n_{X}P=0$ and from the equality
(6.7)(i) we obtain
\[g(AX,Y)\xi + g(Y,\xi)AX=0.\]
We may have one of the following situations:\\
(i) If $AX$ and $\xi$ are linearly dependent vector fields, then
there exist a real number $\A$ such that $AX=\A \xi$ and from this
we obtain $g(Y,\xi)=0$ for any $Y \in \chi(M)$. Thus, for $Y=\xi$
we obtain  $g(\xi,\xi)=0$ which is equivalent with $\xi=0$ and this is an impossible case.\\
(ii) If $AX$ and $\xi$ are linearly independent vector fields,
then $g(AX,Y)=0$ (for any $X,Y \in \chi(M)$). Thus $A=0$ and from
this we have that M is a totally geodesic submanifold in $\MM$. $
\square $

\prop If M is a hypersurface in an almost product Riemannian
manifold $(\MM,\PP,\G)$, with structure $(P,g,\xi,u,a)$ induced
on M by $\PP$, then the following equalities are equivalent:
\[ \n_{X}u=0 \Leftrightarrow \n_{X}\xi=0. \]
\normalfont \textbf{Proof:} If $\n_{X}u=0$ then we obtain
$g(AX,PY)=ag(AX,Y)$, from the equality (6.7)(ii).

From $g(AX,PY)=g(P(AX),Y)$ we have $g(P(AX)-aAX,Y)=0$, for any $X,
Y \in \chi(M)$, and using the equality  (6.7)(iii) we have
$\n_{X}\xi=0$.

Conversely, we suppose that $\n_{X}\xi=0$ (for any $X \in
\chi(M)$) and we have $g(\n_{X}\xi,Y)=0$ from (6.7)(iii). Thus
\[g(P(AX)-aAX, Y)=0 \Longleftrightarrow g(AX,PY)-ag(AX,Y)=0\] for
any $X,Y \in \chi(M)$. Therefore, we obtain $\n_{X}u=0$, from
(ii)(6.7) . $\square $

\normalfont From the theorem 4.2 we have the following property:

\prop Let M be a hypersurface of a locally product Riemannian
manifold $(\MM, \G, \PP)$, with $(P,g,\xi,u,a)$ induced structure
on M by $\PP$. We suppose that $a \neq \pm 1$. A necessary and
sufficient condition for the normality of structure
$(P,g,\xi,u,a)$ is that the tensor field P commutes by the
Weingarten operator A (that is $PA =AP$).

\normalfont From the proposition (3.2) it follows that :

\prop  Let M be a hypersurface in a locally product manifold
$(\MM, \G, \PP)$, with the induced structure $(P,g,\xi,u,a)$. The
necessary and sufficient condition for the 1-form u on M to be
closed (i.e.$du=0$) is that the tensor field P commutes by the
Weingarten operator A.

\prop Let M be a hypersurface in a locally product Riemannian
manifold $(\MM, \G, \PP)$ with $(P,g,\xi,u,a)$ induced structure
on M. Then $\xi$ is a Killing vector field with respect to g on M
if and only if we have
\[2aA=PA+AP \leqno(6.9)\]
where A is the Weingarten operator on M.\\
\normalfont\textbf{Proof:} We have that $\xi$ is a Killing vector
field on M if and only if
\[(L_{\xi}g)(Y,Z)=0, \: (\forall) Y,Z \in \chi(M)\]
\[\Longleftrightarrow \xi g(Y,Z)-g([\xi,Y],Z)-g(Y,[\xi,Z])=0 \]
\[\Longleftrightarrow \xi g(Y,Z)-g(\n_{\xi}Y-\n_{Y}\xi,Z)-g(Y,\n_{\xi}Z-\n_{Z}\xi)=0\]
\[\Longleftrightarrow\underbrace{\xi g(Y,Z)-g(\n_{\xi}Y,Z)-g(Y,\n_{\xi}Z)}_{=0}+g(\n_{Y}\xi,Z)+g(Y,\n_{Z}\xi)=0\]
so
\[g(\n_{Y}\xi,Z)+g(Y,\n_{Z}\xi)=0, \: (\forall) Y,Z \in \chi(M) \leqno(6.10)\]
Using the equality  (6.7)(iii) we obtain that (6.10) is equivalent
with
\[g(-PA_{N}Y+aA_{N}Y,Z)+g(-PA_{N}Z+aA_{N}Z,Y)=0 \]
\[\Longleftrightarrow g(2aA_{N}Y-PA_{N}Y-A_{N}(PY),Z)=0, \: (\forall) Y,Z \in \chi(M)\]
which is equivalent with the equality (6.9). $\square$

\cor Let M be a hypersurface of a locally product Riemannian
manifold $(\MM, \G, \PP)$, with the normal induced structure
$(P,g,\xi,u,a)$. Then $\xi$ is a Killing vector field with respect
to g on M if and only if we have
\[aA=PA=AP \leqno(6.11)\]
where A is the Weingarten operator on M.

\prop Let M be a hypersurface of a locally product Riemannian
manifold $(\MM, \G, \PP)$, with a normal induced structure
$(P,g,\xi,u,a)$ and $\xi$ is a Killing vector field. If $ a^{2}
\neq 1$, then rank A =1 and $\xi$ is
an eigenvector of the Weingarten operator A, with the eigenvalue $\frac{\xi (a)}{2(a^{2}-1)}$.\\
\normalfont \textbf{Proof:} From the corollary 6.1 we have
$PA=aA$. Thus, $P^{2}A=aPA=a^{2}A$ and we have
\[ P^{2}(AX)=a^{2}AX, \: (\forall) X \in \chi(M) \]
Using the equality  (6.3)(i) we obtain  $AX-u(AX)\xi=a^{2}AX$, so
\[ (a^{2}-1)AX=-u(AX)\xi, \: (\forall) X \in \chi(M) \leqno(6.12)\]
From (6.7)(iv) and (6.12) (with $a^{2} \neq 1$) we have
\[ AX=\frac{X(a)}{2(a^{2}-1)}\xi \: (\forall) X \in \chi(M) \leqno(6.13)\]
and from this we remark that rank A=1. More of this, if we put
$X=\xi$ in the equality (6.13) we obtain
\[ A(\xi)=\frac{\xi(a)}{2(a^{2}-1)}\xi \: (\forall) X \in \chi(M) \leqno(6.14)\]
thus $\xi$ is a eigenvector of Weingarten operator A, and its
eigenvalue is $\frac{\xi (a)}{2(a^{2}-1)}$.$\square$

\rem \normalfont Under the assumption of the last proposition, if
\[e_{1}=\frac{\xi}{\sqrt{1-a^{2}}}, \leqno(6.15)\]
then
\[A(e_{1})=A(\frac{\xi}{\sqrt{1-a^{2}}})=\frac{1}{\sqrt{1-a^{2}}}A(\xi)=
\frac{\xi(a)}{2(a^{2}-1)}\frac{\xi}{\sqrt{1-a^{2}}}=\frac{\xi(a)}{2(a^{2}-1)}e_{1}\]
Therefore, $e_{1}$ is eigenvector of the Weingarten operator A and
its eigenvalue is $\frac{\xi(a)}{2(a^{2}-1)}$.

Let $(e_{1},...,e_{n})$ be an orthonormal basis of a tangent space
$T_{x}M$, for any $x\in M$. From (6.13) we have $Ker A =n-1$ and
we obtain $A(e_{i})=0$ for any $i \in \{2,...,n\}$. Thus, the
first eigenvalue of the Weingarten operator is
\[k_{1}=\frac{\xi(a)}{2(a^{2}-1)} \leqno(6.16)\] and the others
are
\[k_{2}=...=k_{n}=0 \leqno(6.17)\]
More of them, the mean vector field H of hipersurface M in $\MM$
is
\[H=\frac{1}{n}trace A \cdot N=\frac{1}{n}k_{1}N\]
and we have
\[H=\frac{\xi(a)}{2n(a^{2}-1)} N. \leqno(6.18)\]

\prop Let M be a hypersurface in a locally product Riemannian
manifold $(\MM, \G, \PP)$, with the normal induced structure
$(P,g,\xi,u,a)$. If we have
$a^{2} \neq 1$ and $\xi(a)=0$, then M is totally geodesic in $\MM$.\\
\normalfont \textbf{Proof:} From (6.16) (with $\xi(a)=0$) and
(6.17) we obtain
\[k_{1}=k_{2}=...=k_{n}=0,\]
so, the Weingarten operator A=0 and we have that M is totally
geodesic in $\MM$. $\square$

\prop Let M be a totally umbilical hypersurface (i.e $A=\lambda
I$) in a locally product  manifold $(\MM, \G, \PP)$, with the
induced structure $(P,g,\xi,u,a)$. Then we have:
\[ \begin{cases}
(i) \quad (\n_{X}P)(Y)=\lambda g(Y,\xi)X+ \lambda g(X,Y)\xi,\\
(ii)\:\:(\n_{X}u)(Y)=-\lambda g(X,PY)+a\lambda g(X,Y),\\
(iii)\:\n_{X}(\xi)=-\lambda PX+a \lambda X, \quad \n_{\xi}\xi=2a \lambda \xi,\\
 (iv)\quad X(a)=-2\lambda g(X,\xi)
 \end{cases}\leqno(6.19)\]
 for any $X,Y \in \chi(M)$.\\
 \normalfont \textbf{Proof:}  We have
 \[g(A_{N}X,Y)=\lambda g(X,Y)\] because M is totally umbilical in
 $\MM$ and using (6.7) we obtain (6.19). $\square$

\prop If M is a totally umbilical hypersurface in a locally
product Riemannian manifold $(\MM, \G, \PP)$, with the induced
structure
$(P,g,\xi,u,a)$, then the 1-form $u$ is closed.\\
 \normalfont \textbf{Proof:} Using the equalities (6.19)(ii) and (3.19), we obtain
\[du(X,Y)=(\n_{X}u)(Y)-(\n_{Y}u)(X)=\lambda(g(PX,Y)-g(X,PY))=0\]
Thus, the 1-form u is closed. $\square$

\cor Let M be a totally umbilical submanifold in a locally product
Riemannian manifold $(\MM, \G, \PP)$, with the induced structure
$(P,g,\xi,u,a)$. Then, it follows that
\[ \begin{cases}
(i) \quad (\n_{X}P)(\xi)= \lambda(1-a^{2})X, \\
(ii) \:\: (\n_{\xi}P)(X)=2\lambda g(X,\xi)\xi,\\
(iii)\:(\n_{X}u)(\xi)=2a \lambda g(X,\xi),
 \end{cases}\leqno(6.20)\]
 for any $X \in \chi(M)$.\\
\normalfont \textbf{Proof:} \\
For $Y=\xi$ in the relations (6.19)(i),(ii) and using the equality
(6.3)(iii), we obtain (i) and (iii) from (6.20). If $X=\xi$ in the
equality  (6.19)(i) we obtain (ii) from (6.20). $\square$

\cor  Let M be a totally umbilical submanifold in a locally
product manifold $(\MM, \G, \PP)$, with the induced structure
$(P,g,\xi,u,a)$ and $X \in \xi^{\bot}$ (where $\xi^{\bot}=\{X \in
\chi(M) / X \bot \xi\}$). Then we have
\[ \begin{cases}
(i) \:\:\: (\n_{\xi}P)(X)=0,\\
(ii)\:\:(\n_{X}u)(\xi)=0,\\
(iii)\: X(a)=0 \: \Longrightarrow \: a=constant
\end{cases}\leqno(6.21)\] for any $X \in \chi(M)$.

\rem \normalfont Let M be a hypersurface in a locally product
Riemannian manifold $(\MM, \G, \PP)$, with the induced structure
$(P,g,\xi,u,a)$. We suppose that $(e_{1},...,e_{n})$ is an
othonormal basis of the tangent space $T_{x}M$, (for any $x \in
M$). Then $div \xi= trace (e_{i} \rightarrow \n_{e_{i}}\xi)$ and
using (6.19)(iii) we obtain $\n_{e_{i}}\xi =
\lambda(aI-P)(e_{i})$. So,
\[div \xi =\lambda (na-trace P) \leqno(6.22)\] Therefore,
if $div \xi =0$ it follows that $trace P =na$.

\bigskip In the following we assume that M is an n-dimensional submanifold of codimension 2 in
an almost product Riemannian manifold $(\MM,\G,\PP)$, with induced
structure $(P,g,u_{\A},\xi_{\A},(a_{\A\B})_{2})$ on M ($\A,\B \in
\{1,2\}$). We suppose that the normal connection vanishes
identically (thus $l_{\A \B}=0$). In these conditions, the
relations (1.6) and (1.7) from Theorem 1.1 have the following
forms:
\[  \begin{cases}
  (i) \quad P^{2}X=X-u_{1}(X)\xi_{1}-u_{2}(X)\xi_{2},  \\
  (ii) \quad u_{1}(PX)=-a_{11}u_{1}(X)-a_{12}u_{2}(X) ,\\
  (iii) \quad  u_{2}(PX)=-a_{21}u_{1}(X)-a_{22}u_{2}(X),\\
  (iv) \quad u_{1}(\xi_{1})=1-a_{11}^{2}-a_{12}^{2},\\
  (v) \quad u_{2}(\xi_{2})=1-a_{12}^{2}-a_{22}^{2},\\
  (vi) \quad  u_{1}(\xi_{2})=u_{2}(\xi_{1})=-a_{12}(a_{11}+a_{22}),\\
  (vii) \quad  P(\xi_{1})= -a_{11}\xi_{1} - a_{12}\xi_{2},\\
  (viii) \quad    P(\xi_{2})= -a_{21}\xi_{1} -a_{22}\xi_{2},\\
  (ix) \quad g(PX,PY)=g(X,Y)-u_{1}(X)u_{1}(Y)-u_{2}(X)u_{2}(Y)
\end{cases}\leqno(6.23)
\]
for any $X, Y \in \chi(M)$.

\normalfont We denote by $\mathcal{A}:=\begin{pmatrix}
  a_{11} & a_{12} \\
  a_{21} & a_{22}
\end{pmatrix}$.

 \rem \normalfont  A simplifier assumption for
these relations is $a_{11}+a_{22}=0$ thus, $trace \:
\mathcal{A}=0$, which is equivalent with $\xi_{1} \bot \xi_{2}$.
Under this assumption, if we denote $a_{11}=-a_{22}=a$,
$a_{12}=a_{21}=b$ and $1-a^{2}-b^{2}=\sigma$, from the relations
(ii)-(vii) (2.1), we easily see that
\[  \begin{cases}
(i) \quad u_{1}(\xi_{1})=u_{2}(\xi_{2})=\sigma \Longleftrightarrow
   g(\xi_{1},\xi_{1})=g(\xi_{2},\xi_{2})=\sigma,\\
(ii) \quad  u_{2}(\xi_{1})=u_{1}(\xi_{2})=0 \Longleftrightarrow
   g(\xi_{1},\xi_{2})=0,\\
(iii) \quad    u_{1}(PX)=-au_{1}(X)-bu_{2}(X),\\
(iv) \quad    u_{2}(PX)=-bu_{1}(X)+au_{2}(X),\\
(v) \quad    P(\xi_{1})= -a\xi_{1} - b\xi_{2},\\
(vi) \quad    P(\xi_{2})= -b\xi_{1} + a\xi_{2}.
  \end{cases} \leqno(6.24)\]

Furthermore, from (2.6), under the assumption that the normal
connection $\n^{\bot}$ vanishes identically (i.e. $l_{\A \B}=0$),
we obtain
\[(\n_{X}P)(Y)=g(A_{1}X,Y)\xi_{1}+g(A_{2}X,Y)\xi_{2}+\leqno(6.25)\]\[+g(Y,\xi_{1})A_{1}X+g(Y,\xi_{2})A_{2}X,\]
\[  \begin{cases}
       (\n_{X}u_{1})(Y)=-g(A_{1}X,PY)+a g(A_{1}X,Y)+b g(A_{2}X,Y),   \\
       (\n_{X}u_{2})(Y)=-g(A_{2}X,PY)+b g(A_{1}X,Y)-a g(A_{2}X,Y) .
  \end{cases}\leqno(6.26)\]and
\[  \begin{cases}
 \n_{X}\xi_{1}=-P(A_{1}X)+a A_{1}X+b A_{2}X,  \\
 \n_{X}\xi_{2}=-P(A_{2}X)+b A_{1}X-a A_{2}X.
 \end{cases}\leqno(6.27) \]and
\[  \begin{cases}
 (i) \quad X(a)=-2g(A_{1}X,\xi_{1}), \\
 (ii) \quad X(b)=-g(A_{1}X,\xi_{2})-g(A_{2}X,\xi_{1})
 \end{cases}\leqno(6.28)\]
for any $X,Y \in \chi(M)$. \\

\lem Let M be an n-dimensional submanifold of codimension 2 in an
 locally product Riemannan manifold $(\MM,\G,\PP)$, with the normal induced structure
 $(P,g,u_{\A},\xi_{\A},(a_{\A\B})_{r})$. If the normal
connection $\n^{\bot}$ vanishes identically (i.e. $l_{\A \B}=0$),
then the following equation is good
\[ g(Y,\xi_{1})B_{1}(X)+g(Y,\xi_{2})B_{2}(X)+
C_{1}(X,Y)\xi_{1}+C_{2}(X,Y)\xi_{2}=0,\leqno(6.29) \] for any $X,Y \in \chi(M)$.\\
\normalfont \textbf{Proof:} By virtue of (3.40) we obtain :
\[ g(X,\xi_{1})B_{1}(Y)+g(X,\xi_{2})B_{2}(Y)=g(Y,\xi_{1})B_{1}(X)+g(Y,\xi_{2})B_{2}(X),\leqno(6.30)\] for any $X,Y
\in \chi(M)$. Multiplying by $Z \in \chi(M)$ the equality (6.30)
we have
\[g(X,\xi_{1})g(B_{1}(Y),Z)+g(X,\xi_{2})g(B_{2}(Y),Z)=\]\[=g(Y,\xi_{1})g(B_{1}(X),Z)+g(Y,\xi_{2})g(B_{2}(X),Z),\] and
using the notation (ii) from (3.36) it follows that
\[
g(X,\xi_{1})C_{1}(Y,Z)+g(X,\xi_{2})C_{2}(Y,Z)=\leqno(6.31)\]\[=g(Y,\xi_{1})C_{1}(X,Z)+g(Y,\xi_{2})C_{2}(X,Z),\]
for any $X,Y,Z \in \chi(M)$. Inverting Y by Z in the last equality
we obtain
\[ g(X,\xi_{1})C_{1}(Z,Y)+g(X,\xi_{2})C_{2}(Z,Y)=\leqno(6.31)'\] \[=g(Z,\xi_{1})C_{1}(X,Y)+g(Z,\xi_{2})C_{2}(X,Y).\]
Summating the relations (6.30) and (6.31)', and using the
skew-symmetry of $C_{1}$ and $C_{2}$ (from (3.36)(ii)), we obtain
\[g(Y,\xi_{1})C_{1}(X,Z)+g(Y,\xi_{2})C_{2}(X,Z)+g(Z,\xi_{1})C_{1}(X,Y)+g(Z,\xi_{2})C_{2}(X,Y)=0\]
which is equivalent with
\[g(Y,\xi_{1})g(B_{1}(X),Z)+g(Y,\xi_{2})g(B_{2}(X),Z)+\]\[+g(Z,\xi_{1})g(B_{1}(X),Y)+g(Z,\xi_{2})g(B_{2}(X),Y)=0\]
Hence,
\[g(g(Y,\xi_{1})B_{1}(X),Z)+g(g(Y,\xi_{2})B_{2}(X),Z)+\]\[+g(Z,g(B_{1}(X),Y)\xi_{1})+g(Z,g(B_{2}(X),Y)\xi_{2})=0\]
and it follows that
\[g([g(Y,\xi_{1})B_{1}(X)+g(Y,\xi_{2})B_{2}(X)+g(B_{1}(X),Y)\xi_{1}+g(B_{2}(X),Y)\xi_{2})],Z)=0\] for any $ Z\in
\chi(M)$ and from this we obtain (6.29). $\square$

\lem Let M be an n-dimensional submanifold of codimension 2 in an
 locally product Riemannan manifold $(\MM,\G,\PP)$, with the normal induced structure
 $(P,g,u_{\A},\xi_{\A},(a_{\A\B})_{r})$. If the normal
connection $\n^{\bot}$ vanishes identically (i.e. $l_{\A \B}=0$)
and $\sigma \neq 0$ then, the following equalities are good
\[ \begin{cases}
  (i) \quad B_{1}(\xi_{1})=0,\\
  (ii)\quad B_{2}(\xi_{2})=0,\\
  (iii)\:\:\: B_{1}(\xi_{2})=0,\\
  (iv) \quad B_{2}(\xi_{1})=0.
  \end{cases} \leqno(6.32)
   \]
\normalfont \textbf{Proof:} By virtue of (6.29) with
$X=Y=\xi_{1}$, we obtain
\[ g(\xi_{1},\xi_{1})B_{1}(\xi_{1})+g(\xi_{1},\xi_{2})B_{2}(\xi_{1})+\leqno(6.33)\]
\[+g(B_{1}(\xi_{1}),\xi_{1})\xi_{1}+g(B_{2}(\xi_{1}),\xi_{1})\xi_{2}=0.\]
Using $g(\xi_{1},\xi_{1})=\sigma \neq 0$, $g(\xi_{1},\xi_{2})=0$,
$g(B_{1}(\xi_{1}),\xi_{1})=C_{1}(\xi_{1},\xi_{1})=0$ and
$g(B_{2}(\xi_{1}),\xi_{1})=C_{2}(\xi_{1},\xi_{1})=0$ we obtain
$B_{1}(\xi_{1})=0$, from the skew-symmetry of $C_{1}$ and $C_{2}$.

From the equality (6.29) with $X=Y=\xi_{2}$, we obtain
\[g(\xi_{2},\xi_{1})B_{1}(\xi_{2})+g(\xi_{2},\xi_{2})B_{2}(\xi_{2})+\leqno(6.34)\]
\[+g(B_{1}(\xi_{2}),\xi_{2})\xi_{1}+g(B_{2}(\xi_{2}),\xi_{2})\xi_{2}=0,\]
Using that $g(\xi_{2},\xi_{2})=\sigma \neq 0$,
$g(\xi_{1},\xi_{2})=0$,
$g(B_{1}(\xi_{2}),\xi_{2})=C_{1}(\xi_{2},\xi_{2})=0$ and
$g(B_{2}(\xi_{2}),\xi_{2})=C_{2}(\xi_{2},\xi_{2})=0$ we obtain
$B_{2}(\xi_{2})=0$, from the skew-symmetry of $C_{1}$ and $C_{2}$.

If we put $X=\xi_{1}$ and $Y=\xi_{2}$ in (6.29), we obtain
\[g(\xi_{2},\xi_{1})B_{1}(\xi_{1})+g(\xi_{2},\xi_{2})B_{2}(\xi_{1})+\leqno(6.35)\]
\[+g(B_{1}(\xi_{1}),\xi_{2})\xi_{1}+g(B_{2}(\xi_{1}),\xi_{2})\xi_{2}=0,\]
Using $g(\xi_{2},\xi_{2})=\sigma \neq 0$, $g(\xi_{1},\xi_{2})=0$,
$B_{1}(\xi_{1})=0$ and
\[g(B_{2}(\xi_{1}),\xi_{2})=C_{2}(\xi_{1},\xi_{2})=
-C_{2}(\xi_{2},\xi_{1})=-g(B_{2}(\xi_{2}),\xi_{1})=0\] and
replacing these in (6.35) we obtain $B_{2}(\xi_{1})=0$.

Using the equality (6.29) with $X=\xi_{2}$ and $Y=\xi_{1}$ we
obtain
\[ g(\xi_{1},\xi_{1})B_{1}(\xi_{2})+g(\xi_{1},\xi_{2})B_{2}(\xi_{2})+\leqno(6.36)\]
\[+g(B_{1}(\xi_{2}),\xi_{1})\xi_{1}+g(B_{2}(\xi_{2}),\xi_{1})\xi_{2}=0,\]
and from $g(\xi_{1},\xi_{1})=\sigma \neq 0$,
$g(\xi_{1},\xi_{2})=0$, $B_{2}(\xi_{2})=0$, $B_{1}(\xi_{1})=0$ and
\[g(B_{1}(\xi_{2}),\xi_{1})=C_{1}(\xi_{2},\xi_{1})=
-C_{1}(\xi_{1},\xi_{2})=-g(B_{1}(\xi_{1}),\xi_{2})=0\] we have
$B_{1}(\xi_{2})=0$. $\square$

 Under the assumption for the codimension r=2, the following proposition can be considered as a particular case
of the theorem 4.1. Using the last two lemmas, we could make
another proof of this, which will be done below:

\prop We suppose that M is a submanifold of codimension 2 in a
 locally product Riemannian manifold $(\MM,\G,\PP)$, with the normal induced structure
 $(P,g,u_{\A},\xi_{\A},(a_{\A\B})_{2})$. If the normal
connection $\n^{\bot}$ vanishes identically (i.e. $l_{\A\B}=0$),
$trace \mathcal{A}=0$ and $\sigma \neq 0$, then P commutes with
the Weingarten operators $A_{\A}$ ($\A \in \{1,2\}$), thus the
following relations take place:
\[  \begin{cases}
       (i) \quad (PA_{1}-A_{1}P)(X)=0,\: (\forall) X \in \chi(M)\\
       (ii) \quad(PA_{2}-A_{2}P)(X)=0,\: (\forall) X \in \chi(M)
\end{cases} \leqno(6.37)\]
\normalfont \textbf{Proof:} With $Y=\xi_{1}$ in the equality
(6.29) we obtain
\[g(\xi_{1},\xi_{1})B_{1}(X)+g(\xi_{1},\xi_{2})B_{2}(X)+\leqno(6.38)\]
\[+g(B_{1}(X),\xi_{1})\xi_{1}+g(B_{2}(X),\xi_{1})\xi_{2}=0,\] and
from $g(\xi_{1},\xi_{1})=\sigma \neq 0$, $g(\xi_{1},\xi_{2})=0$ we
have \[g(B_{\A}(X),\xi_{\B})=C_{\A}(X,\xi_{\B})=
-C_{\A}(\xi_{\B},X)=-g(B_{\A}(\xi_{\B}),X)=0 \leqno(6.39)\] where
$\A, \B \in \{1,2\}$. From the last lemma we have
$B_{\A}(\xi_{\B})=0$, for any $\A, \B \in \{1,2\}$. Therefore we
obtain  $B_{1}X=0$, for any $X \in \chi(M)$, so we have (i) from
(6.37).

With $Y=\xi_{2}$ in (6.29), we obtain
\[g(\xi_{2},\xi_{1})B_{1}(X)+g(\xi_{2},\xi_{2})B_{2}(X)+\leqno(6.40)\]
\[+g(B_{1}(X),\xi_{2})\xi_{1}+g(B_{2}(X),\xi_{2})\xi_{2}=0.\] From
$g(\xi_{2},\xi_{2})=\sigma \neq 0$, $g(\xi_{1},\xi_{2})=0$, using
the equality (6.40) and the relations (ii) and (iii) from the last
lemma then we obtain $B_{2}X=0$, for any $X \in \chi(M)$, so we
have (ii) from (6.37). $ \square $

\prop Let M be a submanifold of codimension 2 in a locally product
Riemannian manifold $(\MM,\G,\PP)$, with the normal induced
structure $(P,g,u_{\A},\xi_{\A},(a_{\A\B})_{2})$. If the normal
connection $\n^{\bot}$ vanishes identically (i.e. $l_{\A\B}=0$),
$trace \mathcal{A}=0$ and $\sigma \neq 0$, then the relations
occur:
\[ \begin{cases}
(i) \quad A_{1}\xi_{1}=\frac{1}{\sigma}h_{1}(\xi_{1},\xi_{1})\xi_{1}+\frac{1}{\sigma}h_{1}(\xi_{1},\xi_{2})\xi_{2},\\
(ii)\quad A_{1}\xi_{2}=\frac{1}{\sigma}h_{1}(\xi_{1},\xi_{2})\xi_{1}+\frac{1}{\sigma}h_{1}(\xi_{2},\xi_{2})\xi_{2},\\
(iii) \quad A_{2}\xi_{1}=\frac{1}{\sigma}h_{2}(\xi_{1},\xi_{1})\xi_{1}+\frac{1}{\sigma}h_{2}(\xi_{1},\xi_{2})\xi_{2},\\
(iv) \quad A_{2}\xi_{2}=\frac{1}{\sigma}h_{2}(\xi_{1},\xi_{2})\xi_{1}+\frac{1}{\sigma}h_{2}(\xi_{2},\xi_{2})\xi_{2}.\\
\end{cases} \leqno(6.41)\]
\normalfont \textbf{Proof:} Applying P in (6.37)(i) it follows
that \[P^{2}A_{1}X=PA_{1}PX, \leqno(6.42)\]  for any $ X \in
\chi(M)$. Using the equality (6.23)(i) we obtain
\[A_{1}X-u_{1}(A_{1}X)\xi_{1}-u_{2}(A_{1}X)\xi_{2}= PA_{1}PX, \leqno(6.43)\]
for any $X \in \chi(M)$.

If we put in (6.43) $X=\xi_{1}$ and $X=\xi_{2}$, respectively,  we
obtain
\[A_{1}\xi_{1}+aPA_{1}\xi_{1}+bPA_{1}\xi_{2}=
h_{1}(\xi_{1},\xi_{1})\xi_{1}+h_{1}(\xi_{1},\xi_{2})\xi_{2}
\leqno(6.44)\] and
\[A_{1}\xi_{2}+bPA_{1}\xi_{1}-aPA_{1}\xi_{2}=
h_{1}(\xi_{1},\xi_{2})\xi_{1}+h_{1}(\xi_{2},\xi_{2})\xi_{2}.
\leqno(6.45)\] from the equalities (6.24)(v) and (vi).

We replace $X\rightarrow PX $ in the equality (6.37)(i) so,
$PA_{1}PX=A_{1}P^{2}X$ and using the equality (6.23)(i), we obtain
\[ PA_{1}PX=A_{1}X-u_{1}(X)A_{1}\xi_{1}-u_{2}(X)A_{1}\xi_{2},
 \leqno(6.46)\] for any $X \in \chi(M)$.

If we put $X=\xi_{1}$ and $X=\xi_{2}$ respectively in (6.46) we
obtain
\[ (\sigma -1)A_{1}\xi_{1}-aPA_{1}\xi_{1}-bPA_{1}\xi_{2}=0, \leqno(6.47)\]
and
\[ (\sigma -1)A_{1}\xi_{2}-bPA_{1}\xi_{1}+aPA_{1}\xi_{2}=0. \leqno(6.48)\]
from the equalities (6.24)(v) and (vi).

Summating the relations (6.44) and (6.47), for $\sigma \neq 0$, we
obtain (i) from (6.41). Summating the relations (6.45) and (6.48)
we obtain, for $\sigma \neq 0$, the equality (ii) from (6.41).

Applying P in the equality (6.37)(ii),it follows that
\[P^{2}A_{2}X=PA_{2}PX,  \leqno(6.49)\] for any $X \in \chi(M)$
and using (6.23)(i), we obtain
\[A_{2}X-u_{1}(A_{2}X)\xi_{1}-u_{2}(A_{2}X)\xi_{2}= PA_{2}PX,  \leqno(6.50)\]
for any $ X \in \chi(M)$.

For $X=\xi_{1}$ and $X=\xi_{2}$, respectively in (6.50), and using
the equalities (6.24)(v) and (vi), we obtain
\[A_{2}\xi_{1}+aPA_{2}\xi_{1}+bPA_{2}\xi_{2}=
h_{2}(\xi_{1},\xi_{1})\xi_{1}+h_{2}(\xi_{1},\xi_{2})\xi_{2}
\leqno(6.51)\] and
\[A_{2}\xi_{2}+bPA_{2}\xi_{1}-aPA_{2}\xi_{2}=
h_{2}(\xi_{1},\xi_{2})\xi_{1}+h_{2}(\xi_{2},\xi_{2})\xi_{2}.
\leqno(6.52)\]

We replace $X\rightarrow PX $ in the equality (6.37)(ii) so,
$PA_{2}PX=A_{2}P^{2}X$ and using the equality (6.23)(i), we obtain
\[ PA_{2}PX=A_{2}X-u_{1}(X)A_{2}\xi_{1}-u_{2}(X)A_{2}\xi_{2},\leqno(6.53)\] for any $X \in \chi(M)$.

For $X=\xi_{1}$ and $X=\xi_{2}$, respectively in (6.53), and using
the equalities (6.24)(v) and (vi), we obtain
\[ (\sigma -1)A_{2}\xi_{1}-aPA_{2}\xi_{1}-bPA_{2}\xi_{2}=0,
\leqno(6.54)\] and
\[ (\sigma -1)A_{2}\xi_{2}-bPA_{2}\xi_{1}+aPA_{2}\xi_{2}=0. \leqno(6.55)\]
Summating  the relations (6.51) and (6.55) we obtain , for $\sigma
\neq 0$, the equality  (iii) from (6.41). Summating the relations
(6.52) and (6.55) we obtain, for $\sigma \neq 0$, the equality
(iv) from (6.41). $\square $


\section{\bf Some examples of structures $(P,g,\xi_{\A},u_{\A},(a_{\A\B})_{r})$ induced  on
 submanifolds of Euclidean  space }
\normalfont

\textbf{ Exemple 1.} We suppose that the ambient space is
$\MM=E^{2p}$ and for any $x \in E^{2p}$ we have
$x=(x^{1},...,x^{p},y^{1},...,y^{p}):=(x^{i},y^{i})$. The tangent
space $T_{x}E^{2p}$ is isomorphic to $E^{2p}$. Let
$\PP:E^{2p}\rightarrow E^{2p}$ an almost product structure on
$E^{2p}$ such that
\[\PP(x^{1},...,x^{p},y^{1},...,y^{p})=(y^{1},...,y^{p},x^{1},...,x^{p}) \leqno(7.2)\]
Thus, $(\PP, <\:>)$ is an almost product Riemannian structure on
$E^{2p}$. We show that any hypersphere $S^{2p-1}\hookrightarrow
E^{2p}$ has an $(a,1)f$ Riemannian structure, by constructing it in a effective way. \\

The equation of sphere $S^{2p-1}(R)$ is
\[(x^{1})^{2}+...+(x^{p})^{2}+(y^{1})^{2}+...+(y^{p})^{2}=R^{2} \leqno(7.3)\]
where R is its radius and $(x^{1},...,x^{p},y^{1},...,y^{p})$ are
the coordinates of any point $x \in S^{2p-1}(R)$. We have that
\[N=\frac{1}{R}(x^{1},...,x^{p},y^{1},...,y^{p}) \leqno(7.4)\]
is a unit normal vector in x on sphere $S^{2p-1}(R)$ and
\[\PP N=\frac{1}{R}(y^{1},...,y^{p},x^{1},...,x^{p}) \leqno(7.5)\]
We denote by $(X^{1},...,X^{p},Y^{1},...,Y^{p})$ a tangent vector
in x at $S^{2p-1}(R)$. Hence we have
\[\sum_{i=1}^{p}x^{i}X^{i}+\sum_{i=1}^{p}y^{i}Y^{i}=0 \leqno(7.6)\]
If we decompose $\PP N$ in the tangential and normal components,
we obtain
\[ \PP N=\frac{1}{R}(\xi^{1},...,\xi^{p},\eta^{1},...,\eta^{p}) + a \cdot
\frac{1}{R}(x^{1},...,x^{p},y^{1},...,y^{p}) \leqno(7.7) \]

From the relations (7.5) and (7.7) we have
\[\xi^{i}=y^{i}-ax^{i}, \quad \eta^{i}=x^{i}-ay^{i} \leqno(7.8)\]
But $(\xi^{1},...,\xi^{p},\eta^{1},...,\eta^{p})$ is tangent at
the sphere $S^{2p-1}(R)$ and from this we obtain
\[\sum_{i=1}^{p}x^{i}\xi^{i}+\sum_{i=1}^{p}y^{i}\eta^{i}=0 \leqno(7.9)\]
Using (7.8), it follows that
\[\sum_{i=1}^{p}x^{i}y^{i}-a\sum_{i=1}^{p}(x^{i})^{2}+
\sum_{i=1}^{p}x^{i}y^{i}-a\sum_{i=1}^{p}(y^{i})^{2}=0\] So, from
the equation (7.3) we have
\[ a=\frac{2}{R^{2}}\sum_{i=1}^{p}x^{i}y^{i} \leqno(7.10)\]
Therefore, the matrix $\mathcal{A}$ becomes a real function $a$ on
$S^{2p-1}(R)$. Moreover, from equality  (7.8)  we obtain the
tangential component of $\PP N$ at sphere $S^{2p-1}(R)$
\[\xi=\frac{1}{R}(y^{1}-ax^{1},...,y^{p}-ax^{p},x^{1}-ay^{1},...,x^{p}-ay^{p}) \leqno(7.11)\]
where $a$ was defined in (7.10).

 Using the equality (7.11) in $u(X)=<X,\xi>$, with
$X=(X^{1},...,X^{p},Y^{1},...,Y^{p})$ a tangent vector in x at
sphere $S^{2p-1}(R)$, then we have
\[u(X)=\frac{1}{R}[\sum_{i=1}^{p}y^{i}X^{i}-a\sum_{i=1}^{p}(x^{i}X^{i}+y^{i}Y^{i})+
\sum_{i=1}^{p}x^{i}Y^{i}]\] and from  (7.6) we obtain
\[ u(X)=\frac{1}{R}(\sum_{i=1}^{p}x^{i}Y^{i}+\sum_{i=1}^{p}y^{i}X^{i}) \leqno(7.12)\]

From $PX=\PP X -u(X)N$ we obtain
\[PX=(Y^{1},...,Y^{p},X^{1},...,X^{p})-
\frac{u(X)}{R}(x^{1},...,x^{p},y^{1},...,y^{p})\] and from this we
have
\[\leqno(7.13) \: PX=(Y^{1}-\frac{u(X)}{R}x^{1},...Y^{p}-\frac{u(X)}{R}x^{p},
X^{1}-\frac{u(X)}{R}y^{1},...,X^{p}-\frac{u(X)}{R}y^{p})\] where
$X=(X^{1},...,X^{p},Y^{1},...,Y^{p})$ is a tangent vector in
$x=(x^{1},...,x^{p},y^{1},...,y^{p})$ at sphere and u(X) was
defined in (7.12). Moreover, $PX$ is tangent at $S^{2p-1}(R)$
because
\[ <PX,N>=\sum_{i=1}^{p}x^{i}Y^{i}-\frac{u(X)}{R}\sum_{i=1}^{p}(x^{i})^{2}+
\sum_{i=1}^{p}y^{i}X^{i}-\frac{u(X)}{R}\sum_{i=1}^{p}(y^{i})^{2}=\]
\[=\sum_{i=1}^{p}(x^{i}Y^{i}+y^{i}X^{i})-\frac{u(X)}{R}\sum_{i=1}^{p}[(x^{i})^{2}+(y^{i})^{2}]=
\sum_{i=1}^{p}(x^{i}Y^{i}+y^{i}X^{i})-u(X)\cdot R =0\] so we have
$<PX,N>=0$.

Furthermore, if X and Y are tangent vectors of the sphere
$S^{2p-1}(R)$, then we have
\[<PX,Y>=<\PP X,Y>=<\PP^{2}X,\PP Y>=<X,\PP Y>=<X,PY>.\]
For $Y=(X'^{1},...,X'^{p},Y'^{1},...,Y'^{p})$ we have $PY =
(Y'^{i}-\frac{u(Y)}{R}x^{i},X'^{i}-\frac{u(Y)}{R}y^{i})$ Thus,
\[<PX,PY>=\sum_{i}(Y^{i}-\frac{u(X)}{R}x^{i})(Y'^{i}-\frac{u(Y)}{R}x^{i})+
\sum_{i}(X^{i}-\frac{u(X)}{R}y^{i})(X'^{i}-\frac{u(Y)}{R}y^{i})=\]
\[=\sum_{i}(X^{i}X'^{i}+Y^{i}Y'^{i})+\frac{u(X)u(Y)}{R^{2}}\underbrace{
\sum_{i}((x^{i})^{2}+(y^{i})^{2})}_{R^{2}}-\]
\[-\frac{u(X)}{R}\underbrace{\sum_{i}(x^{i}Y'^{i}+y^{i}X'^{i})}_{R\cdot u(Y)}-
\frac{u(Y)}{R}\underbrace{\sum_{i}(x^{i}Y^{i}+y^{i}X^{i})}_{R\cdot
u(X)}\] and from this we have $<PX,PY>=<X,Y>-u(X)u(Y)$, for any
tangent vectors X and Y on the sphere $S^{2p-1}(R)$, in any point
$x \in M$.

Therefore, from the relations (7.10), (7.11), (7.12), (7.13) we
have a $(P,\xi,u,a)$ induced structure by $\PP$ from $E^{2p}$ on
the sphere $S^{2p-1}(R)$. This structure is a normal $(a,1)f$
Riemannian structure, because  $S^{2p-1}(R)$ is the totally
umbilical hypersurface in $E^{2p}$ and from this we have that the
tensor field P commutes with the Weingarten operator A.

\textbf{ Example 2.} Let $S^{2p-1}(1)$ a hypersphere of the
Euclidean  space $E^{2p}$ ($p \geq 2$), endowed with an almost
product Riemannian structure $(\PP, <\:>)$ from the previous
example. We have seen, from above, that on any hypersphere
$S^{2p-1}(1)$ we have an induced structure $(P,\xi,u,a)$. It is
obvious that $E^{2p}=E^{p} \times E^{p}$ and in each of $E^{p}$ we
can get a hypersphere
\[S^{p-1}(r_{1})=\{(x^{1},...,x^{p}), \: \sum_{i=1}^{p}(x^{i})^{2}=r_{1}^{2}\}, \leqno(7.14)\]
and
\[S^{p-1}(r_{2})=\{(y^{1},...,y^{p}), \: \sum_{i=1}^{p}(y^{i})^{2}=r_{2}^{2}\},\leqno(7.15)\]
respectively, with the assumption that
$r_{1}^{2}+r_{2}^{2}=1$.

We construct the product manifold $S^{p-1}(r_{1}) \times
S^{p-1}(r_{2})$ (such in \cite{IS}, pg.116-117).

Any point $x \in S^{p-1}(r_{1}) \times S^{p-1}(r_{2})$ has the
coordinates $(x^{1},...,x^{p},y^{1},...,y^{p})$ with the
properties $\sum_{i=1}^{p}(x^{i})^{2}=r_{1}^{2}$ and
$\sum_{i=1}^{p}(y^{i})^{2}=r_{2}^{2}$. This manifold is a
submanifold of codimension 2 in $E^{2p}$. Furthermore,
$S^{p-1}(r_{1}) \times S^{p-1}(r_{2})$ (with
$r_{1}^{2}+r_{2}^{2}=1$) is a submanifold of codimension 1 in
$S^{2p-1}(1)$. Therefore we have the successive imbedded
\[S^{p-1}(r_{1}) \times S^{p-1}(r_{2}) \hookrightarrow S^{2p-1}(1) \hookrightarrow E^{2p} \leqno(7.16)\]

The tangent space in a point $x=(x^{1},...,x^{p},y^{1},...,y^{p})$
at the product of spheres $S^{p-1}(r_{1}) \times S^{p-1}(r_{2})$
is
\[T_{(x^{1},...,x^{p},0,...,0)} S^{p-1}(r_{1}) \oplus T_{(0,...,0,y^{1},...,y^{p})}S^{p-1}(r_{2})\]

A vector $(U^{1},...,U^{p})$ from $T_{x}E^{p}$ (with
$x=(x^{1},...,x^{p})$) is tangent to $S^{p-1}(r_{1})$ if
\[\sum_{i=1}^{p}U^{i}x^{i}=0 \leqno(7.17)\]
and it can be identified by $(U^{1},...,U^{p},0,...,0)$ from
$E^{2p}$.

A vector $(V^{1},...,V^{p})$ from $T_{y}E^{p}$ (with
$y=(y^{1},...,y^{p})$) is tangent to $S^{p-1}(r_{2})$ if
\[\sum_{i=1}^{p}V^{i}y^{i}=0 \leqno(7.18)\]
and it can be identified by $(0,...,0,V^{1},...,V^{p})$ from
$E^{2p}$.

Consequently, for any $(x,y):=(x^{1},...,x^{p},y^{1},...,y^{p})
\in S^{p-1}(r_{1}) \times S^{p-1}(r_{2})$ we have
\[(U^{1},...,U^{p},V^{1},...,V^{p}):=(U^{i},V^{i}) \in T_{(x,y)}(S^{p-1}(r_{1}) \times
S^{p-1}(r_{2})) \] if and only if
\[\sum_{i=1}^{p}U^{i}x^{i}=0, \quad \sum_{i=1}^{p}V^{i}y^{i}=0.\leqno(7.19)\]
Furthermore, from (7.19) we remark that
$\sum_{i=1}^{p}(U^{i}x^{i}+V^{i}y^{i})=0$ so $(U^{i},V^{i})$ is a
tangent vector at $S^{2p-1}(1)$ in any $(x,y)\in S^{p-1}(r_{1})
\times S^{p-1}(r_{2})$. From this it follows that
\[T_{(x,y)}(S^{p-1}(r_{1}) \times S^{p-1}(r_{2})) \subset T_{(x,y)}S^{2p-1},\leqno(7.20)\]
for any $(x,y) \in S^{p-1}(r_{1}) \times S^{p-1}(r_{2})$.

We denote by $N_{2}$ the normal unit vector at $S^{2p-1}(1)$ in a
point (x,y). Thus, we have
\[N_{2}=(x^{1},...,x^{p},y^{1},...,y^{p}), \leqno(7.21)\]
from
$\sum_{i=1}^{p}(x^{1})^{2}+\sum_{i=1}^{p}(y^{1})^{2}=r_{1}^{2}+r_{2}^{2}=1$.
On the other hand, $N_{2}$ is a normal vector field at
$(S^{p-1}(r_{1}) \times S^{p-1}(r_{2}))$, when it is considered in
its points. Let
\[N_{1}=(\frac{r_{2}}{r_{1}}x^{1},...,\frac{r_{2}}{r_{1}}x^{p},
-\frac{r_{1}}{r_{2}}y_{1},...,-\frac{r_{1}}{r_{2}}y_{p}):=
(\frac{r_{2}}{r_{1}}x^{i},-\frac{r_{1}}{r_{2}}y_{i}),\leqno(7.22)\]
be a vector in any $(x,y) \in S^{p-1}(r_{1}) \times
S^{p-1}(r_{2})$. We remark that
\[<N_{1},N_{2}> = \frac{r_{2}}{r_{1}}\sum_{i}(x^{i})^{2}-\frac{r_{1}}{r_{2}}\sum_{i}(y^{i})^{2}=
\frac{r_{2}}{r_{1}}r_{1}^{2}-\frac{r_{1}}{r_{2}}r_{2}^{2}=0,\]
Therefore $N_{2}$ is orthogonal to $N_{1}$ in any $(x,y) \in
S^{p-1}(r_{1}) \times S^{p-1}(r_{2})$ and from this we have that
$N_{1}$ is a tangent vector at $S^{2p-1}(1)$. Moreover, $N_{1}$ is
a unit vector because
\[<N_{1},N_{1}>=\frac{r_{2}^{2}}{r_{1}^{2}}
\sum_{i}(x^{i})^{2}+\frac{r_{1}^{2}}{r_{2}^{2}}\sum_{i}(y^{i})^{2}=r^{2}_{2}+r_{1}^{2}=1\]
in any $(x,y) \in S^{p-1}(r_{1}) \times S^{p-1}(r_{2})$.

Let $U=(U^{1},...,U^{p},V^{1},...,V^{p}):=(U^{i},V^{i})$ be a
tangent vector at $S^{p-1}(r_{1}) \times S^{p-1}(r_{2})$ in any
its point $(x^{i}, y^{i})$. From (7.19) we have
\[<U,N_{1}>= \frac{r_{2}}{r_{1}}\sum_{i}x^{i}U^{i}-\frac{r_{1}}{r_{2}}\sum_{i}y^{i}V^{i}=0\]
so $N_{1}$ is a normal vector field at $S^{p-1}(r_{1}) \times
S^{p-1}(r_{2})$ in any $(x,y) \in S^{p-1}(r_{1}) \times
S^{p-1}(r_{2})$, and $(N_{1},N_{2})$ is an orthonormal basis in
$T_{(x^{i},y^{i})}^{\bot}S^{p-1}(r_{1}) \times S^{p-1}(r_{2})$ in
any point $(x,y) \in S^{p-1}(r_{1}) \times S^{p-1}(r_{2})$.

We have the induced structure $(P,\xi,u,a)$ on $S^{2p-1}(1)$ which
was constructed in (7.10), (7.11), (7.12), (7.13). In the
following we find the induced structure on $S^{p-1}(r_{1}) \times
S^{p-1}(r_{2})$ by the structure $(P,\xi,u,a)$ on $S^{2p-1}(1)$,
using the propositions 5.1 and 5.2. Thus, we shall have a
$(P_{0},\xi_{0},\xi^{\top},u_{0},u,(a_{\A\B}))$ induced structure
on the submanifold $S^{p-1}(r_{1})\times S^{p-1}(r_{2})$ by the
$(P,\xi,u,a)$ structure on $S^{2p-1}(1)$.

Using the relations (7.12) and (7.22) we have
\[u(N_{1})=\sum_{i}(-x^{i}y^{i}\frac{r_{1}}{r_{2}}+\frac{r_{2}}{r_{1}}x^{i}y^{i})=
(\frac{r_{2}}{r_{1}}-\frac{r_{1}}{r_{2}})\sum_{i}x^{i}y^{i}\] We
denote by $\lambda:=\frac{r_{2}}{r_{1}}-\frac{r_{1}}{r_{2}}$ and
by $\sigma:=\sum_{i}x^{i}y^{i}$. Then, it follows that
\[u(N_{1})=\lambda \sigma \leqno(7.23)\]

From the equality (7.10), with $R=1$, we have
\[a=2\sum_{i=1}^{p}x^{i}y^{i}=2\sigma. \leqno(7.24)\]

If we decomposed $P(N_{1})$ in normal and tangential components at
$S^{p-1}(r_{1}) \times S^{p-1}(r_{2})$ in $S^{2p-1}(1)$ we obtain
\[P(N_{1})=(\xi^{1},...,\xi^{p},\eta^{1},...,\eta^{p})+b N_{1} \leqno(7.25)\]
where $(\xi^{1},...,\xi^{p},\eta^{1},...,\eta^{p})$ is a tangent
field at $S^{p-1}(r_{1}) \times S^{p-1}(r_{2})$ and $b$ is a real
function on this submanifold. Using the equality (7.13) we obtain
\[P(N_{1})=(-\frac{r_{1}}{r_{2}}y^{1}-\lambda \sigma x^{1},...,-\frac{r_{1}}{r_{2}}y^{p}-\lambda \sigma
x^{p}, \frac{r_{2}}{r_{1}}x^{1}-\lambda \sigma
y^{1},...,\frac{r_{2}}{r_{1}}x^{p}-\lambda \sigma y^{p} )
\leqno(7.26)\]
 Thus, from (7.22), (7.25) and (7.26), we obtain
\[ \begin{cases}
(i) \:\: \xi^{i}=-b\frac{r_{2}}{r_{1}}x^{i}-\frac{r_{1}}{r_{2}}y^{i}-\lambda \sigma x^{i},\\
 (ii)\: \eta^{i}=b\frac{r_{1}}{r_{2}}y^{i}+\frac{r_{2}}{r_{1}}x^{i}-\lambda \sigma
 y^{i}. \end{cases}\leqno(7.27)\]
 Hence, from (7.27)(i) and using that $\sum_{i}\xi^{i}x^{i}=0$
(because  $(\xi^{1},...,\xi^{p},\eta^{1},...,\eta^{p})$ is tangent
to $S^{p-1}(r_{1}) \times S^{p-1}(r_{2})$), we obtain
 \[b r_{1}r_{2}=-\frac{r_{1}}{r_{2}}\sigma - (\frac{r_{2}}{r_{1}}-\frac{r_{1}}{r_{2}})r_{1}^{2}\sigma\]
Furthermore, from $r_{1}^{2}+r_{2}^{2}=1$ we have
 \[b=-2\sigma. \leqno(7.28)\]

From (7.27) and (7.28), the tangential component at $P(N_{1})$ is
\[ \xi_{0}=(\xi^{1},...,\xi^{p},\eta^{1},...,\eta^{p}):= (\xi^{i},\eta^{i})=(\frac{\sigma}{r_{1}r_{2}}x^{i}-\frac{r_{1}}{r_{2}}y^{i},
\frac{r_{2}}{r_{1}}x^{i}-\frac{\sigma}{r_{1}r_{2}}y^{i}).\leqno(7.29)\]

From $u_{0}(X)=<X,\xi_{0}>$ (where $X=(X^{i},Y^{i})$ is a tangent
vector field of $S^{p-1}(r_{1}) \times S^{p-1}(r_{2})$), we can
find the 1-form $u_{0}$. So, using $\sum_{i}X^{i}x^{i}=0$ and
$\sum_{i}Y^{i}y^{i}=0$, we obtain
\[ u_{0}(X)=\frac{r_{2}}{r_{1}}\sum_{i}x^{i}Y^{i}-\frac{r_{1}}{r_{2}}\sum_{i}y^{i}X^{i}.\leqno(7.30)\]

We denoted by $P_{0}$ the tangent component of the (1,1) tensor
field P (defined in (7.13)) at $S^{p-1}(r_{1}) \times
S^{p-1}(r_{2})$. For the tangent vector field
$(X^{1},...,X^{p},Y^{1},...,Y^{p}):=(X^{i},Y^{i})$ at
$S^{p-1}(r_{1}) \times S^{p-1}(r_{2})$, we have
\[P_{0}(X^{i},Y^{i})=P(X^{i},Y^{i})-u_{0}(X^{i},Y^{i})N_{1} \leqno(7.31)\]
From (7.13), (7.22) and (7.31) it follows that
\[ P_{0}(X^{i},Y^{i})=(Y^{i}-\frac{1}{r_{1}^{2}}(\sum_{j=1}^{p}x^{j}Y^{j})x^{i},
X^{i}-\frac{1}{r_{2}^{2}}(\sum_{j=1}^{p}X^{j}y^{j})y^{i})
\leqno(7.32)\]

On the other hand, from (5.11) we have
\[\xi^{\bot}=<\xi,N_{1}>N_{1}=u(N_{1})N_{1}\]
and from (7.23) it follows that
\[\xi^{\bot}=\lambda \sigma N_{1}= \sigma
((\frac{r_{2}^{2}}{r_{1}^{2}}-1)x^{i},(\frac{r_{1}^{2}}{r_{2}^{2}}-1)y^{i})\leqno(7.33)\]
 From $\xi^{\top}=\xi-\xi^{\bot}$ and using the relations (7.11) and (7.33) we obtain
\[\xi^{\top}=(y^{i}-\frac{\sigma}{r_{1}^{2}}x^{i},x^{i}-\frac{\sigma}{r_{2}^{2}}y^{i}) \leqno(7.34)\]
From (5.10) we can find the entries $a_{12}=a_{21}$ of the $2
\times 2$ matrix $\mathcal{A}$. Hence,
\[a_{12}=<\xi,N_{1}>=u(N_{1})=\lambda \sigma \leqno(7.35)\]
 Therefore, from the relations (7.24), (7.28) and (7.35)
 we obtain the matrix $\mathcal{A}=(a_{\A\B})$ with its entries
\[a_{11}=a=2\sigma, \quad a_{22}=b=-2\sigma, \quad a_{12}=a_{21}=\lambda \sigma\]
 which is
\[\mathcal{A}=\begin{pmatrix}
  2\sigma &  \lambda \sigma\\
  \lambda \sigma & -2\sigma
\end{pmatrix}\leqno(7.36)\]
where $\lambda=\frac{r_{2}}{r_{1}}-\frac{r_{1}}{r_{2}}$ and
$\sigma=\sum_{i}x^{i}y^{i}$.

Consequently, from the corollary 5.1 we obtain the
$(P_{0},\xi_{0},\xi^{\top},u_{0},u,(a_{\A\B}))$ induced structure
on $S^{p-1}(r_{1}) \times S^{p-1}(r_{1})$ by the almost product
Riemannian structure $(\PP, <>)$ on $E^{2p}$, which is effectively
determined by the relations (7.11),( 7.12), (7.30), (7.32),
(7.33), (7.36) and it is an $(a,1)f$ Riemannian structure.

\bigskip

\textbf{ Example 3.} We suppose that the ambient space is the
Euclidean  space $E^{2p+1}$ and let $\PP$ be an almost product
Riemannian structure defined by
\[\PP:E^{2p+1}\rightarrow E^{2p+1}\]
\[\PP(x^{1},...,x^{p},t,y^{1},...,y^{p})=(y^{1},...,y^{p},t,x^{1},...,x^{p}) \leqno(7.37)\]
We show that any hypersphere $S^{2p}(R)$ in $E^{2p+1}$ has a
normal $(a,1)f$ Riemannian structure. The equation of sphere
$S^{2p}(R)$ in
$x=(x^{1},...,x^{p},t,y^{1},...,y^{p}):=(x^{i},y^{i})$ is
\[(x^{1})^{2}+...+(x^{p})^{2}+t^{2}+(y^{1})^{2}+...+(y^{p})^{2}=R^{2} \leqno(7.38)\]
where R is its radius. Then, we remark that
\[N=\frac{1}{R}(x^{1},...,x^{p},t,y^{1},...,y^{p}) \leqno(7.39)\]
is a unit normal vector on $S^{2p}(R)$ in any point $x \in
S^{2p}(R)$. From this we have
\[\PP N=\frac{1}{R}(y^{1},...,y^{p},t,x^{1},...,x^{p}) \leqno(7.40)\]
Let $(X^{1},...,X^{p},T,Y^{1},...,Y^{p})$ be a tangent vector
field on  $S^{2p}(R)$. Thus, from (7.39) it follows that
\[\sum_{i=1}^{p}x^{i}X^{i}+tT+\sum_{i=1}^{p}y^{i}Y^{i}=0. \leqno(7.41)\]
We decomposed $\PP N$ in the tangential and normal components:
\[ \PP N=\frac{1}{R}(\xi^{1},...,\xi^{p},\tau,\eta^{1},...,\eta^{p}) + a \cdot
\frac{1}{R}(x^{1},...,x^{p},t,y^{1},...,y^{p}) \leqno(7.42) \]

From the relations (7.40) and (7.42) we obtain
$y^{i}=\xi^{i}+ax^{i}$, $t=\tau+ a\cdot t$ and
$x^{i}=\eta^{i}+ay^{i}$ (for $i \in \{1,...,p\}$), so
\[\xi^{i}=y^{i}-ax^{i}, \quad \tau =t(1-a), \quad \eta^{i}=x^{i}-ay^{i} \leqno(7.43)\]
But $(\xi^{1},...,\xi^{p},\tau,\eta^{1},...,\eta^{p})$ must to be
tangent on $S^{2p}$, thus it follows that
\[\sum_{i=1}^{p}x^{i}\xi^{i}+\tau \cdot t+\sum_{i=1}^{p}y^{i}\eta^{i}=0 \leqno(7.44)\]
and from this, we obtain
\[\sum_{i=1}^{p}x^{i}y^{i}-a\sum_{i=1}^{p}(x^{i})^{2}+t^{2}(1-a)+
\sum_{i=1}^{p}x^{i}y^{i}-a\sum_{i=1}^{p}(y^{i})^{2}=0\] and from
this we have
\[ a=\frac{1}{R^{2}}(2\sum_{i=1}^{p}x^{i}y^{i}+t^{2}) \leqno(7.45)\]

From (7.43), we have
\[\xi=\frac{1}{R}(y^{1}-ax^{1},...,y^{p}-ax^{p},t(1-a),
x^{1}-ay^{1},...,x^{p}-ay^{p}) \leqno(7.46)\] where $a$ was
defined in (7.45).

Using (7.46) in $u(X)=<X,\xi>$, where
$X=(X^{1},...,X^{p},T,Y^{1},...,Y^{p})$ is a tangent vector field
on sphere, we have
\[u(X)=\frac{1}{R}[(\sum_{i=1}^{p}y^{i}X^{i}+\sum_{i=1}^{p}x^{i}Y^{i}+Tt)-
a\underbrace{(\sum_{i=1}^{p}(x^{i}X^{i}+y^{i}Y^{i})+tT)}_{=0}] \]
and from (7.41), we obtain
\[ u(X)=\frac{1}{R}(\sum_{i=1}^{p}x^{i}Y^{i}+\sum_{i=1}^{p}y_{i}X^{i}+tT) \leqno(7.47)\]

From $PX=\PP X -u(X)N $ we obtain
\[PX=(Y^{1},...,Y^{p},T,X^{1},...,X^{p})-
\frac{u(X)}{R}(x^{1},...,x^{p},t,y^{1},...,y^{p})\] Therefore, it
follows that
\[PX=(Y^{i}-\frac{u(X)}{R}x^{i},T-\frac{u(X)}{R}t,X^{i}-\frac{u(X)}{R}y^{i}), \leqno(7.48)\]
for $i\in \{1,...,p\}$.

We verify that $PX$ is tangent at $S^{2p-1}$. Using the relations
(7.48), we obtain
\[ <PX,N>=\sum_{i=1}^{p}x^{i}Y^{i}-\frac{u(X)}{R}\sum_{i=1}^{p}(x^{i})^{2}+(T-\frac{u(X)}{R}t)t+
\sum_{i=1}^{p}y^{i}X^{i}-\frac{u(X)}{R}\sum_{i=1}^{p}(y^{i})^{2}=\]
\[=\sum_{i=1}^{p}(x^{i}Y^{i}+tT+y^{i}X^{i})-
\frac{u(X)}{R}\sum_{i=1}^{p}[(x^{i})^{2}+t^{2}+(y^{i})^{2}]=\]
\[=\sum_{i=1}^{p}(x^{i}Y^{i}+tT+y^{i}X^{i})-u(X)\cdot R =0\]
thus $ <PX,N>=0$ so $PX$ is tangent sphere $S^{2p}$.

Furthermore, if X and Y are tangent vector fields on sphere
$S^{2p}$, then
\[<PX,Y>=<\PP X,Y>=<\PP^{2}X,\PP Y>=<X,\PP Y>=<X,PY>\]
If $Y:=(X'^{1},...,X'^{p},T',Y'^{1},...,y'^{p})$ is a tangent
vector in any point x at sphere $S^{2p}$ then, from (7.50), we
have
\[PY = (Y'^{i}-\frac{u(Y)}{R}x^{i},T'-\frac{u(Y)}{R}t,X'^{i}-\frac{u(Y)}{R}y^{i}),\:\:
i \in \{1,...,p\} \leqno(7.49)\] and from (7.48) and (7.49) we
have
\[<PX,PY>=\sum_{i}(Y^{i}-\frac{u(X)}{R}x^{i})(Y'^{i}-\frac{u(Y)}{R}x^{i})+\]
\[+(T-\frac{u(X)}{R}t)(T'-\frac{u(Y)}{R}t)+\sum_{i}(X^{i}-\frac{u(X)}{R}y^{i})(X'^{i}-\frac{u(Y)}{R}y^{i})=\]
\[=\sum_{i}(X^{i}X'^{i}+TT'+Y^{i}Y'^{i})+\frac{u(X)u(Y)}{R^{2}}\underbrace{
\sum_{i}((x^{i})^{2}+t^{2}+(y^{i})^{2})}_{R^{2}}-\]
\[-\frac{u(X)}{R}\underbrace{\sum_{i}(x^{i}Y'^{i}+tT'+y^{i}X'^{i})}_{R\cdot u(Y)}-
\frac{u(Y)}{R}\underbrace{\sum_{i}(x^{i}Y^{i}+tT+y^{i}X^{i})}_{R\cdot
u(X)}\] Hence, we have
\[<PX,PY>=<X,Y>-u(X)u(Y)\]
for any tangent vector fields  X and Y on $S^{2p-1}$.

Consequently, from the relations (7.45), (7.46), (7.47), (7.48) we
have that $(P,\xi,u,a)$ is an induced structure on $S^{2p}$ by the
almost product Riemannian structure ($\PP,<>$) from $E^{2p+1}$.
This structure induced on $S^{2p}$ is a normal $(a,1)f$-Riemannian
structure, because the sphere $S^{2p}$ is a totally umbilical
hypersurface in $E^{2p+1}$ thus, P commutes by the Weingarten
operator A.

\bigskip
\textbf{ Example 4.} Let $\MM =E^{p+q}$ be the ambient space and
we define an almost product Riemannian structure $(\PP,<>)$  on
the Euclidean  space $E^{p+q}$ such that $\PP:E^{p+q}\rightarrow
E^{p+q}$ and
\[\PP(x^{1},...,x^{p},y^{1},...,y^{q})=(x^{1},...,x^{p},-y^{1},...,-y^{q}) \leqno(7.50)\]
We show that any hypersphere $S^{p+q-1}(R)$ in $E^{p+q}$ has an
$(a,1)f$ Riemannian induced structure by $\PP$.

We denote by $(x^{1},...,x^{p},y^{1},...,y^{q}):=(x^{i},y^{j})$,
with $i \in\{1,...,p\}$ and $j \in \{1,...,q\}.$

The equation of the sphere $S^{p+q-1}(R)$ in a point
$x=(x^{i},y^{j})$ is
\[(x^{1})^{2}+...+(x^{p})^{2}+(y^{1})^{2}+...+(y^{q})^{2}=R^{2} \leqno(7.51)\]

The  unit normal vector at sphere $S^{p+q-1}(R)$ is
\[N=\frac{1}{R}(x^{1},...,x^{p},y^{1},...,y^{q}) \leqno(7.52)\]
for any point $x=(x^{i},y^{j})\in S^{p+q-1}$. From (7.50) and
(7.52) we have
\[\PP N=\frac{1}{R}(x^{1},...,x^{p},-y^{1},...,-y^{q}) \leqno(7.53)\]
We denoted by $(X^{1},...,X^{p},Y^{1},...,Y^{q})$ a tangent vector
field at $S^{p+q-1}$. This is orthogonal on N, so we obtain
 \[\sum_{i=1}^{p}x^{i}X^{i}+\sum_{j=1}^{q}y^{i}Y^{i}=0 \leqno(7.54)\]

We decompose $\PP N$ in the tangential and normal components:
\[ \PP N=\underbrace{\frac{1}{R}(\xi^{1},...,\xi^{p},\eta^{1},...,\eta^{q})}_{\xi} + a \cdot
\underbrace{\frac{1}{R}(x^{1},...,x^{p},y^{1},...,y^{q})}_{N}
\leqno(7.55) \]

From (7.52) and (7.54) we obtain
\[\xi^{i}=(1-a)x^{i}, \quad \eta^{j}=-(1+a)y^{j} \leqno(7.56)\]
for any $i \in\{1,...,p\}$ and $j \in \{1,...,q\}$.

Because $(\xi^{1},...,\xi^{p},\eta^{1},...,\eta^{q})$ is tangent
at $S^{p+q-1}$, it follows that
\[\sum_{i=1}^{p}x^{i}\xi^{i}+\sum_{j=1}^{q}y^{j}\eta^{j}=0 \leqno(7.57)\]
From (7.56), we have
\[(1-a)\sum_{i=1}^{p}(x^{i})^{2}-(1+a)\sum_{j=1}^{q}(y^{j})^{2}=0\]
and using the equality (7.51) we obtain
\[ a=\frac{1}{R^{2}}[\sum_{i=1}^{p}(x^{i})^{2}-\sum_{j=1}^{q}(y^{j})^{2}] \leqno(7.58)\]
Therefore, the matrix $\mathcal{A}$, is reduces to a real function
on $S^{p+q-1}$.

From (7.55) we have
\[\xi=\frac{1}{R}((1-a)x^{i},-(1+a)y^{j}) \leqno(7.59)\]
for $i \in\{1,...,p\}$ and $j \in \{1,...,q\}.$

Using (7.59) in $u(X)=<X,\xi>$, with
$X=(X^{1},...,X^{p},Y^{1},...,Y^{q})$ a tangent vector field of
sphere, we have
\[u(X)=\frac{1}{R}[(1-a)\sum_{i=1}^{p}x^{i}X^{i}-(1+a)
\sum_{j=1}^{q}y^{j}Y^{j}] \] and from  (7.54) we obtain
\[ u(X)=\frac{1}{R}(\sum_{i=1}^{p}x^{i}X^{i}-\sum_{j=1}^{q}y^{j}Y^{j}) \leqno(7.60)\]

Let $X=(X^{1},...,X^{p},Y^{1},...,Y^{p})$ a tangent vector field
on sphere $S^{p+q-1}(R)$. From $PX=\PP X -u(X)N$ we have
\[PX=(X^{1},...,X^{p},-Y^{1},...,-Y^{q})-
\frac{u(X)}{R}(x^{1},...,x^{p},y^{1},...,y^{p})\] and from this we
obtain
\[PX=(X^{i}-\frac{u(X)}{R}x^{i},-Y^{j}-\frac{u(X)}{R}y^{j})\leqno(7.61)\]
for $i \in \{1,...,p\}$ and $j \in \{1,...,q\}$. Furthermore, $PX$
is tangent at $S^{p+q-1}$, because we have
\[<PX,N>=(\sum_{i=1}^{p}x^{i}X^{i}-\sum_{j=1}^{q}y^{j}Y^{j})-
\frac{u(X)}{R}[\sum_{i=1}^{p}(x^{i})^{2}+\sum_{j=1}^{q}(y^{j})^{2}]=0\]
On the other hand, if X and Y are tangent vector fields on sphere
$S^{p+q-1}$, then
\[<PX,Y>=<\PP X,Y>=<\PP^{2}X,\PP Y>=<X,\PP Y>=<X,PY>\]
Let $Y=(X'^{1},...,X'^{p},Y'^{1},...,Y'^{q})$ be a tangent vector
in a point $x=(x^{i},y^{j})$ at sphere. From (7.61) we have
\[PY = (X'^{i}-\frac{u(Y)}{R}x^{i},-Y'^{j}-\frac{u(Y)}{R}y^{j})\]
for $ i \in \{1,...,p\}$ and $j \in \{1,...,q\}$. Thus, we have
\[<PX,PY>=\sum_{i=1}^{p}(X^{i}-\frac{u(X)}{R}x^{i})(X'^{i}-\frac{u(Y)}{R}x^{i})+
\sum_{j=1}^{q}(-Y^{j}-\frac{u(X)}{R}y^{j})(-Y'^{j}-\frac{u(Y)}{R}y^{j})=\]
\[=(\sum_{i}X^{i}X'^{i}+\sum_{j=1}^{q}Y^{j}Y'^{j})+
\frac{u(X)u(Y)}{R^{2}}\underbrace{
[\sum_{i=1}^{p}((x^{i})^{2}+\sum_{j=1}^{q}(y^{j})^{2})}_{R^{2}}-\]
\[-\frac{u(X)}{R}\underbrace{(\sum_{i=1}^{p}x^{i}X'^{i}-\sum_{j=1}^{q}y^{j}Y'^{j})}_{R\cdot u(Y)}-
\frac{u(Y)}{R}\underbrace{(\sum_{i}^{p}x^{i}X^{i}-\sum_{j=1}^{q}y^{j}Y^{j})}_{R\cdot
u(X)}\] so
\[<PX,PY>=<X,Y>-u(X)u(Y),\] for any tangent vector fields X and Y on sphere
$S^{p+q-1}$.

Consequently, the relations (7.58), (7.59), (7.60), (7.61) give
the structure $(P,\xi,u,a)$ induced on the sphere $S^{p+q-1}(R)$
by $\PP$ from $E^{p+q}$. This structure is a normal $(a,1)f$
Riemannian structure because, the sphere $S^{p+q-1}$ is a
hypersurface totally umbilical in $E^{p+q}$ thus, P commutes by
the Weingarten operator A.


\bigskip

 "\c{S}tefan cel Mare" University, Suceava, Romania\\
  E-mail address: cristinah@usv.ro


\begin{thebibliography}{n}

\bibitem{Adati1}
T. Adati, T. Miyazawa, Hypersurfaces immersed in an almost
product Riemannian manifold, TRUMath., 14-2 (1978), 17-26.

\bibitem{Adati2}
T. Adati, Submanifolds of an almost product Riemannian manifold,
Kodai Math. J., 4-2 (1981), 327-343.

\bibitem{Adati3}
T. Adati, T. Miyazawa, On paracontact Riemannian manifolds, TRU
Math., 13-2 (1977), 27-39.

\bibitem{AM1}
M. Anastasiei, Some Riemannian almost product structures on
tangent manifold, Proceedings of the 11th National Conference on
Finsler, Lagrange and Hamilton Geometry (Craiova, 2000). Algebras
Groups Geom. 17 (2000), no. 3, 253-262.

\bibitem{Atceken}
M. At\c{c}eken, S. Kele\c{s}, On Product Riemannian Manifolds,
Differential Geometry-Dynamical Systems, Vol.5, No.1, 2003, 1-7.

\bibitem{Atceken2}
M. At\c{c}eken, S. Kele\c{s}, Two Theorems On Invariant
Submanifolds of Riemannian Product Manifold, Indian J. Pure and
Appl. Math., 34(7),July 2003, 1035-1044.

\bibitem{Blair1}
D.E.Blair, G.D. Ludden, K. Yano,Induced structures on
submanifolds, Kodai Math. Sem. Rep. 22 (1970), 188-198.

\bibitem{Blair2}
D.E.Blair, G.D Ludden, Hypersurfaces in almost contact manifolds,
Thoku Math.J.22 (1969),354-362.

\bibitem{Bucki1}
A.Bucki, A. Miernovski, Almost r-paracontact structure, Ann. Univ.
Marie Curie Sklodowska, Sect. A, 39 (1985), 13-26.

\bibitem{Bucki2}
A. Bucki, A. Miernovski, Invariant hypersurfaces of an lmost
r-paracontact manifold, Demonstr. Math., 19(1), (1986), 113-121.

\bibitem{BYC}
B.Y. Chen, Geometry of Submanifolds, Marcel Dekker Inc, New York
(1973).

\bibitem{Chern1}
S.S. Chern, J.G. Wolfson, Minimal surfaces by moving frames, Amer.
J. Math. 105 (1983).

\bibitem{Cruceanu}
V. Cruceanu, Almost Hyperproduct structures on manifolds, An.
Stiint. Univ. Al. I. Cuza Ia»si. Mat. (N. S.) 48 (2002), no. 2,
337-354.

\bibitem{Gray}
A. Gray, Pseudo-Riemannian almost product manifolds and
submersions, J. Math. Mech. 16 (1967), 715-737.


\bibitem{eu5}
C.E. Hre\c{t}canu, On submanifolds in riemannian almost product
manifolds, MASEE International Congress on Mathematics MICOM 2006,
Cyprus, MASEE International Congress on Mathematics - Abstracts,
ISBN 9963-9277-0-X)

\bibitem{eu6}
C.E. Hre\c{t}canu, A structure of paracontact type on submanifolds
in almost product manifolds, 5 th Symposium on Mathematics Applied
in Biology and Biophysics, U.A.S.V.M Iasi, 16-17 june 2006

\bibitem{IS}
S. Ianu\c{s}, Geometrie diferen\c{t}ial\u{a} cu aplica\c{t}ii
\^{i}n teoria relativit\u{a}\c{t}ii, Editura Academiei RSR,
Bucure\c{s}ti 1983.

\bibitem{IS2}
S. Ianu\c{s}, Some almost product structures on manifolds with
linear connection, Kodai Math. Sem. Rep.  23, no. 3 (1971),
305–310.

\bibitem{IS3}
S. Ianu\c{s}, Some submanifolds of tangent bundles, The Procced.
Nat. Sem. on Finsler Spaces, 1980, p.69-76.

\bibitem{IS7}
S. Ianu\c{s}, I.Mihai, On the semi invariant submanifolds in an
almost paracontact riemannian manifolds, Homagial Vol. Kawaguchi,
Tensor 39(1982), p.195-200.

\bibitem{Matsumoto}
K. Matsumoto, On Submanifolds of Locally Product Riemannian
Manifolds, TRU Mathematics 18-2, 1982, 145-157.

\bibitem{IMihai4}
I. Mihai, C.Nicolau, Almost product structures on the tangent
bundle of an almost paracontact manifold, Demonstratio Math.,
15(1982), 1045-1058.

\bibitem{IMihai7}
I. Mihai, C.Nicolau, Some structures induced on the tangent bundle
of an almost paracontact manifold, Proc. Nat. Sem. Finsler Spaces,
Brasov, 1983, 141-146.

\bibitem{IMihai9}
I. Mihai,S.Ianus, K.Matsumoto, Almost semi-invariant submanifolds
of some almost paracontact Riemannian manifolds, Bull. Yamagata
Univ.,11, 1985, 121-128.

\bibitem{IMihai12}
I. Mihai, R. Ro\c{s}ca, L. Verstraelen, Some aspects of the
differential geometry of vector fields, Katholicke Universiteit
Leuven, PADGE, vol 2, 1996.

\bibitem{Miyazawa}
T. Miyazawa, Hypersurfaces immersed in an almost product
Riemannian manifold, Tensor N. S., 33-1 (1979), 114-116.


\bibitem{Naveira}
A.N. Naveira, A classification of Riemannian almost product
structures, Rend. Mat Roma 3 (1983), 577-592.


\bibitem{Nikic}
J. Nikic, Conditions for invariant submanifold of a manifold with
the $(\phi,\xi,\eta,G)$-structure, Kragujevac J. Math. 25 (2003),
pag 147-154.

\bibitem{Okumura}
M. Okumura, Totally umbilical hypersurfaces of a locally product
Riemannian manifold, Kodai Math. Sem. Rep., 19 (1967), 35-42.

\bibitem{Pitis1}
G. Piti\c{s}, On some submanifolds of a locally product manifold,
Kodai Math.J. 9 (1986), 327—333

\bibitem{Sato}
I. Sato, On a structure similar to the almost contact structure
I;II, Tensor N.S. 30 (1976), 219-224; 31 (1977), 199-205.

\bibitem{Sahin}
B. \c{S}ahin, M. At\c{c}eken, Semi-Invariant Submanifolds of
Riemannian Product Manifold, Balkan Journal of Geometry and Its
Applications, Vol.8, No.1, 2003, pp. 91-100.

\bibitem{Senlin}
X. Senlin and N. Yilong, Submanifolds of Product Riemannian
Manifold, Acta Mathematica Scientia 2000, 20(B) 213-218.

\bibitem{Tachibana}
S. Tachibana, Some theorems on locally product Riemannian spaces,
Thoku Math. Jour., 12 (1960), 281-292.

\bibitem{Walker}
A. G. Walker, Almost-product structures, Proc. Sympos. Pure Math.
3 (1961), 94–100.

\bibitem{KY}
K. Yano, M. Kon, Structures on Manifolds, World Scientific,
Singapore, Series in pure matematics - Volume 3, 1984.

\bibitem{Yano1}
K. Yano, M. Okumura, On normal (f, g, u, v, $\lambda$)-structures
on submanifolds of codimension 2 in an even-dimensional Euclidean
space,  Kodai Math. Sem. Rep 23 (1971), 172-197.

\bibitem{Yano2}
K. Yano, M. Okumura, On (f, g, u, v, $\lambda$) -structures, Kodai
Math. Sem. Rep. 22 (1970), 401-423.



\end{thebibliography}
\end{document}